\pdfoutput=1
\RequirePackage{ifpdf}
\ifpdf 
\documentclass[pdftex]{sigma}
\else
\documentclass{sigma}
\fi

\numberwithin{equation}{section}

\usepackage[mathscr]{euscript}
\usepackage[all]{xy}
\usepackage{ytableau}
\usepackage{enumitem}

\newcommand{\C}{\mathbb{C}}

\newcommand{\co}{\xi}

\newcommand{\GL}{\mathrm{GL}}

\newcommand{\Gr}{\mathrm{Gr}}
\newcommand{\Gro}[1]{\mathfrak{G}^Q_{#1}}

\newcommand{\Iaf}{\tilde{I}}

\newcommand{\la}{\lambda}

\newcommand{\loc}{\mathrm{loc}}

\renewcommand{\O}{\mathscr{O}}

\newcommand{\PGL}{\mathrm{PGL}}

\newcommand{\reg}{\Delta}
\newcommand{\Rep}{{R(T)}}

\newcommand{\SL}{\mathrm{SL}}

\newcommand{\ks}{g}
\newcommand{\kst}{\tilde{g}}

\newcommand{\wz}{w_\circ}
\newcommand{\Z}{\mathbb{Z}}
\newcommand{\de}{\varrho}
\newcommand{\sinf}{\frac{\infty}{2}}
\newcommand{\fac}{\psi}

\newcommand{\ha}{f}
\usepackage{chngcntr}

\def\CC{\mathord{\mathbb{C}}}

\def\ZZ{\mathord{\mathbb{Z}}}

\newcommand{\bm}[1]{\mbox{\boldmath{$#1$}}}

\newtheorem{Theorem}{Theorem}[section]
\newtheorem*{Theorem*}{Theorem}
\newtheorem{Corollary}[Theorem]{Corollary}
\newtheorem{Lemma}[Theorem]{Lemma}
\newtheorem{Proposition}[Theorem]{Proposition}

\theoremstyle{definition}
\newtheorem{Definition}[Theorem]{Definition}

\newtheorem{Example}[Theorem]{Example}
\newtheorem{Remark}[Theorem]{Remark}

\begin{document}

\allowdisplaybreaks

\newcommand{\arXivNumber}{2505.02941}

\renewcommand{\PaperNumber}{065}

\FirstPageHeading

\ShortArticleName{Relativistic Toda Lattice and Equivariant $K$-Homology of Affine Grassmannian}

\ArticleName{Relativistic Toda Lattice and Equivariant\\ $\boldsymbol{K}$-Homology of Affine Grassmannian}

\Author{Takeshi IKEDA~$^{\rm a}$, Shinsuke IWAO~$^{\rm b}$, Satoshi NAITO~$^{\rm c}$ and Kohei YAMAGUCHI~$^{\rm a}$}

\AuthorNameForHeading{T.~Ikeda, S.~Iwao, S.~Naito and K.~Yamaguchi}

\Address{$^{\rm a)}$~Faculty of Science and Engineering, Waseda University,\\
\hphantom{$^{\rm a)}$}~3-4-1 Okubo, Shinjuku-ku, Tokyo 169-8555, Japan}
\EmailD{\mail{gakuikeda@waseda.jp}, \mail{yamaguchi\_86@aoni.waseda.jp}}

\Address{$^{\rm b)}$~Faculty of Business and Commerce, Keio University,\\
\hphantom{$^{\rm b)}$}~4-1-1 Hiyosi, Kohoku-ku, Yokohama-si, Kanagawa 223-8521, Japan}
\EmailD{\mail{iwao-s@keio.jp}}

\Address{$^{\rm c)}$~Department of Mathematics, Institute of Science Tokyo,\\
\hphantom{$^{\rm c)}$}~2-12-1 Oh-Okayama, Meguro-ku, Tokyo 152-8551, Japan}
\EmailD{\mail{naito@math.titech.ac.jp}}

\ArticleDates{Received October 26, 2025, in final form June 18, 2026; Published online July 06, 2026}

\Abstract{We investigate the phenomenon known as ``quantum equals affine'' in the setting of~$T$-equivariant quantum $K$-theory of the flag variety $G/B$, as established by Kato for any semisimple algebraic group~$G$. In particular, we focus on the $K$-Peterson isomorphism between the $T$-equivariant quantum $K$-ring $QK_T(\mathrm{SL}_n(\mathbb{C})/B)$ and the $T$-equivariant $K$-homology ring $K_*^T(\mathrm{Gr}_{\mathrm{SL}_n})$ of the affine Grassmannian, after suitable localizations on both sides. Building on an earlier work by Ikeda, Iwao, and Maeno, we present an explicit algebraic realization of the $K$-Peterson map via a~rational substitution that sends the generators of the quantum $K$-theory ring to explicit rational expressions in the fundamental generators of $K_*^T(\mathrm{Gr}_{\mathrm{SL}_n})$, thereby matching the Schubert bases on both sides. Our approach builds on recent developments in the theory of $QK_T(\mathrm{SL}_n(\mathbb{C})/B)$ by Maeno, Naito, and Sagaki, as well as the theory of $K$-theoretic double $k$-Schur functions introduced by Ikeda, Shimozono, and Yamaguchi. This concrete formulation provides new insight into the combinatorial structure of the $K$-Peterson isomorphism in the equivariant setting. As an application, we establish a~factorization formula for the $K$-theoretic double $k$-Schur function associated with the maximal $k$-irreducible $k$-bounded partition.}


\Classification{14N15; 05E10; 37K10}

\vspace{-2mm}

\section{Introduction}

Let $G$ be a~simple, simply-connected algebraic group over $\C$. Fix a~Borel subgroup $B$ of~$G$, and a~maximal torus $T$ contained in $B$.
We study a~remarkable relation
between the $T$\nobreakdash-equivariant
 quantum $K$-ring $QK_T(G/B)$
and
the $T$-equivariant $K$-homology ring
$K_*^T(\Gr_G)$ of the affine Grassmannian $\Gr_G$ of $G$.
This viewpoint, commonly referred to as the ``quantum equals affine'' phenomenon, was originally introduced by
Peterson~\cite{Pet} in the (co)homology context (see also~\cite{LS:Acta}),
and its $K$-theory analogue
has since been investigated by several authors (see the next paragraph for references).

There exists a~map, known as the $K$-\emph{Peterson map}, which connects the rings
$QK_T(G/B)$ and~$K_*^T(\Gr_G)$, after appropriate localization.
The purpose of this paper is
to study the $K$\nobreakdash-Pe\-ter\-son map for $G=\SL_n$ by realizing it through an explicit rational substitution, which establishes a~correspondence between the Schubert bases
at the combinatorial level.

A heuristic version of ``quantum equals affine'' phenomena in $K$-theory was explored for $G=\SL_n$
in non-equivariant setting~\cite{IIM}. This approach utilized
an integrable system, called the \emph{relativictic Toda lattice}, due to Ruijsenaars~\cite{Rui}, which
can be seen as the group version of the ordinary
Lie algebra version of the Toda lattice.
The construction
in~\cite{IIM}
was somewhat \emph{ad hoc} as it relied on an unsolved conjecture by
Kirillov and Maeno concerning the ring presentation of
$QK(\SL_n(\C)/B)$ at that time.
However, this conjecture has since been resolved, with a~minor correction, by Maeno, Naito, and Sagaki~\cite{MNS1}; we also refer the reader to~\cite{ACT,HK,Kim,KPSZ} (see Remark~\ref{rem:KPSZ} below). This advancement has improved the situation significantly.
The $K$-theoretic Peterson isomorphism
for general $G$ was conjectured by Lam, Li, Mihalcea, and Shimozono~\cite{LLMS}, and later proved by Kato~\cite{Kato} by using semi-infinite flag manifolds.
An alternative proof was also provided by Chow and Leung~\cite{CL}.
It is shown in~\cite{IIN} that the map studied in~\cite{IIM} coincides with the
map by Kato, up to a~natural ring automorphism
$\sigma$ (see Section~\ref{sec:sigma}).

The aim of this article is to extend the study in~\cite{IIM} to the equivariant case, building on more recent developments. Specifically,
Maeno, Naito, and Sagaki~\cite{MNS1,MNS2} established that the quantum double Grothendieck polynomials $\Gro{w}$ of Lenart--Maeno represent the Schubert classes $\mathscr{O}^w\in QK_T(\SL_n(\C)/B)$.
On the affine side, Ikeda, Shimozono, and Yamaguchi~\cite{ISY} provided a~realization of $K_*^T(\Gr_{\SL_n})$ in terms of equivariantly deformed symmetric functions.
They introduced a~family of special functions, \smash{$\tilde{g}_x^{(k)}(y|b)$}, called the \emph{$K$-theoretic double $k$-Schur functions}, which are identified with the Schubert classes $\mathscr{O}_x$.
The primary goal of this paper is to realize the $K$-Peterson map explicitly through an algebraic substitution.
This provides a~concrete connection between these Schubert representatives in both the quantum and affine settings.

\subsection[K-theoretic Peterson map: abstract form]{$\boldsymbol{K}$-theoretic Peterson map: Abstract form}
To describe the $K$-Peterson map more precisely, we fix some notation of the affine Weyl group.
Let $W_G$ be the Weyl group of $(G,T)$, and $\hat{W}_G=W\ltimes Q^\vee$ the affine Weyl group,
 where $Q^\vee$ is the coroot lattice.
 Let $\hat{W}_G^0$ be the set of minimal-length coset representatives for
 $\hat{W}_G/W_G$. For~\smash{$x\in \hat{W}_G^0$}, there is
 an associated Schubert structure sheaf $\mathscr{O}_x$ in \smash{$K_*^T(\Gr_{G})$}.
 These sheaves form an $\Rep$-basis of $K_*^T(\Gr_{G})$, where $\Rep$ denotes the representation ring of $T$.
On the quantum side, we have the Schubert class $\mathscr{O}^w$ for each $w\in W_G$.
 For \smash{$x\in \hat{W}_G^0$},
 write it as~${x=wt_\xi}$ with~${w\in W}$, $\xi\in Q^\vee$,
 where $t_\xi$ is the translation element corresponding to $\xi\in Q^\vee$.
 Let $QK_T(G/B)_Q$ denote the localization of $QK_T(G/B)$ by
 the Novikov variables $Q_1,\dots,Q_r$, where~$r$ is the rank of $G$.
 The $K$-Peterson map, at the abstract level of Schubert bases,
maps $Q^\xi\mathscr{O}^w \in QK_T(G/B)_Q$ to $\mathscr{O}_{x}\in K_*^T(\Gr_G)$, where $Q^\xi$ denote the product of the Novikov variables corresponding to $\xi$. This map establishes a~correspondence between the quantum and affine Schubert calculus.

\subsection[K-Peterson map from the Relativistic Toda lattice]{$\boldsymbol{K}$-Peterson map from the Relativistic Toda lattice}

Let us explain the key idea of our construction before going into the details.
By solving the relativistic Toda lattice, we obtain a~birational map between the phase space $\mathscr{Y}$ of the relativistic Toda lattice
and a~certain centralizer family $\mathscr{Z}$ associated with $\PGL_n(\C)$.
Both varieties are defined over $T$, the maximal torus of $\SL_n(\C)$, and their coordinate rings have the following geometric interpretations:
\begin{equation}\label{eq:O_geom}
\mathscr{O}(\mathscr{Y})
\cong
QK_T^{\mathrm{pol}}(\SL_n(\C)/B),\qquad
\mathscr{O}(\mathscr{Z})
\cong
K_*^T(\Gr_{\SL_n}).
\end{equation}
Here $QK_T^{\mathrm{pol}}(\SL_n(\C)/B)$ is the polynomial version of the
quantum $K$-ring of the flag manifold (see Section~\ref{sec:QK}).
The second isomorphism of~\eqref{eq:O_geom} was proved in~\cite{ISY}, which is a~$K$-theory analogue of a~result for $H_*^T(\Gr_G)$ due to Peterson~\cite{Pet} and Ginzburg~\cite{Gi} independently.
We obtain the following commutative diagram (see Theorem~\ref{thm:main})
\begin{equation}
\label{eq:comdiag}
\begin{split}&
\xymatrix{
QK_T^{\mathrm{pol}}(\SL_n(\C)/B)_{Q}
\ar[d] \ar[r] \ar@{}[dr] & K_*^T(\Gr_{\SL_n})\big[{\sigma_i^{-1},\tau_i^{-1}\mid 1\leq i\leq n}\big]
 \ar[d] \\
\mathscr{O}(\mathscr{Y}^\circ) \ar[r]_{\widetilde{\Phi}_n} & {\mathscr{O}}(\mathscr{Z}^\circ),
}\end{split}\end{equation}
 where the top arrow is Kato's map, and the bottom arrow is defined by~\eqref{eq:tPhi} below.
In this diagram, all maps are isomorphisms, and $\mathscr{Y}^\circ$ and $\mathscr{Z}^\circ$ are certain open dense subsets of $\mathscr{Y}$ and~$\mathscr{Z}$ respectively, and $\sigma_i$, $\tau_i$ are discussed in detail in Section~\ref{sec:tau}.
This perspective extends the work~\cite{LS:double-Kostant,LS10} by Lam and Shimozono for
(co)homology and the classical Toda lattice.

\subsection[K-theoretic double k-Schur function]{$\boldsymbol{K}$-theoretic double $\boldsymbol{k}$-Schur function}
\label{sec:double_K-k}

Another ingredient of
our work is the \emph{$K$-theoretic double $k$-Schur function} introduced in~\cite{ISY}, which is an equivariant deformation of the $K$-theoretic $k$-Schur function introduced by Lam, Schilling, and Shimozono~\cite{LSS:K}.
The representation ring $\Rep$ is given as
\[
\Z\big[{\rm e}^{\pm a_i} (1\le i\le n)\big]/\bigl({\rm e}^{a_1+\cdots+a_n}-1\bigr).
\]
For \smash{$x\in \hat{W}_G^0$}, the $K$-theoretic double $k$-Schur function \smash{$g_x^{(k)}(y|b)$} is a~symmetric formal power series in the infinitely many variables $y=(y_1,y_2,\dots)$ with coefficients in~$\Rep$, where we set $k=n-1$.
It depends on the sequence $b=(b_1,\dots,b_n)$ of equivariant parameters, where each~$b_i$~is identified
with $1-{\rm e}^{-a_i}\in \Rep$.
 Let \smash{$\hat{\Lambda}^{\Rep}_{(n)}$} be the $\Rep$-span of \smash{$g_x^{(k)}(y|b)$}, $
 x\in \hat{W}_G^0$. Then we have an isomorphism
\[
 K_*^T(\Gr_{\SL_n})\cong \hat{\Lambda}^{\Rep}_{(n)}
\]
 of $\Rep$-algebras such that the structure sheaf~$\mathscr{O}_x$ corresponds to the \emph{closed} $K$-theoretic double $k$-Schur function~\smash{$\tilde{g}_x^{(k)}(y|b):=\sum_{z\le x}g_z^{(k)}(y|b)$}, where $\le $ denotes
 the Bruhat order on $\hat{W}_G^0$.

There is a~bijection $\hat{W}_G^0\cong \mathscr{P}^{(k)}$, where
$\mathscr{P}^{(k)}$ denotes the set of $k$-bounded partitions, i.e., the~partition $\la=(\la_1,\dots,\la_i)$ such that
$\la_1\le k$.
If \smash{$x\in \hat{W}_G^0$} corresponds to
\smash{$\la\in \mathscr{P}^{(k)}$},
then we write \smash{$g_x^{(k)}(y|b)$} \big(resp.\ \smash{$\tilde{g}_x^{(k)}(y|b)$}\big)
as \smash{$g_\la^{(k)}(y|b)$} \big(resp.\ \smash{$\tilde{g}_\la^{(k)}(y|b)$}\big).

We have derived determinantal formulas for
\smash{$g^{(k)}_\la(y|b)$} and \smash{$\tilde{g}_\la^{(k)}(y|b)$} (see Theorems~\ref{thm:determinant_formula_for_ksmall_g} and~\ref{thm:determinantal_formula_for_g})
for a~$k$-\emph{small} $k$-bounded partition $\la$; a~partition $\la\in \mathscr{P}^{(k)}$ is said to be $k$-small if $\la_1+\ell(\la)\le n$, where
$\ell(\la)$ is the number of nonzero parts of $\la$.

\subsection[k-rectangles and the tau-functions]{$\boldsymbol{k}$-rectangles and the $\boldsymbol{\tau}$-functions}
The so-called $\tau$-\emph{functions} in the theory of integrable systems also play an important role in our geometric context. For $1\le i\le n$, we define
$\tau_i,\sigma_i\in \mathscr{O}(\mathscr{Z})$ as the $i$-th principal minor determinants of certain matrices
related to the centralizer family (see Section~\ref{sec:tau}).

For $1\le i\le k$, let $R_i$ denote the partition $(i)^{n-i}\in \mathscr{P}^{(k)}$ of rectangular shape, and $R_n=\varnothing$.
Note that each $R_i$ $(1\le i\le n-1)$ is a~maximal $k$-small $k$-bounded partition.
Under the isomorphism
\smash{$ \mathscr{O}(\mathscr{Z})
\cong \hat{\Lambda}^{\Rep}_{(n)}$},
$\tau_i$ and $\sigma_i$
correspond to \smash{$g_{R_{n-i}}^{(k)}(y|b)$} and \smash{$\tilde{g}_{R_{n-i}}^{(k)}(y|b)$} up to
some simple factors respectively (see Corollary~\ref{cor:k rect and tau}).

These functions are fundamental because
\smash{$\tilde{g}_{\la}^{(k)}(y|b)$} satisfies the $k$-rectangle factorization property (see Theorem~\ref{thm:$k$-rect})
\[
\tilde{g}_{R_i\cup \la}^{(k)}(y|b)=\tilde{g}_{R_i}^{(k)}(y|b)
\tilde{g}_{\la}^{(k)}\bigl(y|\omega^ib\bigr),
\]
where $\omega$ is the permutation sending $b_i$ to $b_{i+1}$ with $b_{n+1}=b_1$.
Thanks to this result, we only need to study the functions \smash{$\tilde{g}_\la^{(k)}(y|b)$} associated with the
$k$-\emph{irreducible} $k$-bounded partition $\la$; i.e., an element of $\mathscr{P}^{(k)}$ which is not expressed as $R_i\cup\mu$ for $1\le i\le n-1$ and $\mu\in \mathscr{P}^{(k)}$
with~${\mu\ne \varnothing}$.

\subsection{Quantum double Grothendieck polynomials}

For $w\in S_n$, let $\Gro{w}(z| \eta)$ be the \emph{quantum double
Grothendieck polynomial} due to Lenart and Maeno~\cite{LM}.
This is a~polynomial in two sets of variables
$z_1,\dots,z_n$ and $\eta_1,\dots,\eta_n$ with coefficients in
$\Z[Q_1,\dots,Q_{n-1}]$.
We basically follow the notation in~\cite{MNS2}, however,
there are some differences in the identification of the equivariant parameters.
In particular, $\eta_i$ is identified with $1-{\rm e}^{a_{n-i+1}}\in \Rep$ (see~\eqref{eq:eta-b} for more details).
We denote $\Rep[z_1,\dots,z_n,Q_1,\dots,Q_{n-1}]$ by $\Rep[z,Q]$.

In the context of the relativistic Toda lattice,
$z_i$, $Q_i$ are interpreted as dynamical variables (see Section~\ref{sec:Rel Toda}).
There are conserved quantities $F_i(z,Q)\in \Rep[z,Q]$ $(1\le i\le n)$ of the relativistic Toda lattice.
Let $\mathscr{I}_n^{Q,\mathrm{pol}}$ be the
ideal generated by $F_i(z,Q)-e_i({\rm e}^{-a_1},\dots,{\rm e}^{-a_n})$ $(1\le i\le n)$, where $e_i$ denotes the $i$-th elementary
symmetric polynomial.
The ring \smash{$\Rep[z,Q]/\mathscr{I}_n^{Q,\mathrm{pol}}$} is by our definition
the ring of regular functions $\mathscr{O}(\mathscr{Y})$ on the
phase space $\mathscr{Y}$.

Let $\mathscr{Y}^\circ$ be the open set of $\mathscr{Y}$ defined as the complement of the divisor given by the equation $Q_1\cdots Q_{n-1}=0$.
Due to~\cite{MNS2}, the ring $\mathscr{O}(\mathscr{Y}^\circ)=\mathscr{O}(\mathscr{Y})\big[Q_i^{-1}\mid 1\leq i\leq n-1\big]$
is identified with~\smash{$QK_T^{\mathrm{pol}}(\SL_n(\C)/B)_Q$}
and \smash{$\mathfrak{G}_w^Q(z| \eta)$} represents the Schubert structure sheaf \smash{$\mathscr{O}^w$} (see Remark~\ref{rem:QK=O(Y)} for more details).

\subsection{Correspondence of Schubert bases}
Following an analogous approach to that used in Kostant's construction of solutions to the Toda lattice, we obtain a~map
\[\Phi_n\colon \ \mathscr{O}(\mathscr{Y}^\circ)
\longrightarrow \hat{\Lambda}^{\Rep}_{(n)}\big[\tau_i^{-1},\sigma_i^{-1}\mid 1\leq i\leq n\big]
\cong\mathscr{O}(\mathscr{Z}^\circ)\]
of $\Rep$-algebras as
\begin{equation}
z_i\mapsto \frac{\tau_i\sigma_{i-1}}{\sigma_i \tau_{i-1}},\qquad
Q_i\mapsto \frac{\tau_{i-1}\tau_{i+1}}{\tau_i^2}.\label{eq:Phi}
\end{equation}

Let $\omega_k$, with $k=n-1$, be an involution on $\hat{W}_G^0$ called the $k$-\emph{conjugation} defined by replacing~$s_i$ with $s_{n-i}$ for $i \in I:=\{1,\dots,n-1\}$ in any reduced expression of \smash{$x\in \hat{W}_G^0$}.
There is an automorphism $\sigma$ of \smash{$\hat{\Lambda}^{\Rep}_{(n)}$} as an~$\Rep$-algebra
sending $h_i(y)\ (i\in\mathbb{N})$ to $1+h_1(y)+\cdots+h_i(y)$, where~$h_i(y)$ is the
$i$-th complete symmetric function.
In particular, we have $\sigma(\tau_i)=\sigma_i$.
Define the map $\tilde{\Phi}_n$ as follows
\begin{equation}
\label{eq:tPhi}
\tilde{\Phi}_n=\sigma \circ \Phi_n.
\end{equation}
The main result of this paper is the following.
\begin{Theorem}\label{thm:main}
We have the commutative diagram~\eqref{eq:comdiag}. More precisely,
for \smash{$x\in \hat{W}_G^0$}, write $x=wt_\xi$ with $w\in W$, $\xi\in Q^\vee$. Then
\[
\tilde{\Phi}_n\bigl(
Q^\xi \Gro{w}(z| \eta)\bigr)=
\tilde{g}_{x^{\omega_k}}^{(k)}(y|b).
\]
\end{Theorem}
Thus, the substitution map ${\Phi}_n$ agrees with Kato's map up to a~twist, after identifying
 $QK_T(\SL_n(\C)/B)$ with $\mathscr{O}(\mathscr{Y})$ \big(resp.\ $K_*^T(\Gr_{\SL_n})$ with $\O(\mathscr{Z})$\big).

Let us outline the proof.
We begin by proving the theorem for $x=s_0$ in a~purely combinatorial manner (see Proposition~\ref{prop:key}).
In this proof, we utilize basic properties of the quantum double Grothendieck polynomials, the realization of $K_*^T(\Gr_{SL_n})$ as $\mathscr{O}(\mathscr{Z})$, and the explicit description of $\Phi_n$.
This step forms the technical core of the paper.
An essential idea to prove the general case is the use of the action of Demazure operators $D_x^Q$, $x\in \hat{W}_G^0$, on $QK_T(\SL_n(\C)/B)_Q$.

As a~straightforward consequence of the theorem for the specific case
$x=s_0$,
we demonstrate that the map $\tilde{\Phi}_n$
intertwines the actions of the Demazure operators on both the quantum and affine sides (see Corollary~\ref{cor:Phi_commutes_w_Dem}).
The final result needed to complete the proof is the
precise formula for the action of
$D_i^Q$ on the Schubert structure sheaves $\mathscr{O}^w$ in
$QK_T(\SL_n(\C)/B)$
(see Proposition~\ref{prop:left_Dem}).

\subsection[Factorization formula for the maximal k-irreducible partition nu\_n]{Factorization formula for the maximal $\boldsymbol{k}$-irreducible partition $\boldsymbol{\nu_n}$}\label{ssec:fac max}
As an application of
Theorem~\ref{thm:main}, we derive a~factorization formula of \smash{$\tilde
{g}_{\nu_n}^{(k)}(y|b)$}, where
\[
\nu_n=\bigcup_{i=1}^{n-2}(n-i-1)^i.
\]
This partition $\nu_n$ is important because it is the unique maximal $k$-irreducible $k$-bounded partition.
The non-equivariant version of the following result was conjectured in~\cite{IIM}, and it was proved by Blasiak, Morse, and Seelinger~\cite{BMS}.

Let $\lfloor x\rfloor$ denote the greatest integer less than or equal to $x\in \mathbb{R}$.
\begin{Theorem}\label{thm:max factor}
Let $n$ be even and write $n = 2m$. Then
\begin{align}
\label{eq:even}
\tilde{g}^{(k)}_{\nu_n}(y|b)
= \prod_{i=0}^{\lfloor (m-1)/2 \rfloor}
 \frac{\Omega(b_{m-1-2i} | y)}{\Omega(b_{m+2-2i} | y)}
 \prod_{i=1}^{n-2}
 \tilde{g}^{(k)}_{(n-i-1)^i}\bigl(y | \omega^{m+2i+1} b\bigr),
\end{align}
where
$\Omega(b_i | y)$ is defined in~\eqref{eq:s0_Omega}.
Let $n$ be odd. Then
\begin{align}
\label{eq:odd}
\tilde{g}^{(k)}_{\nu_n}(y|b)
= \prod_{i=1}^{n-2}
 \tilde{g}^{(k)}_{(n-i-1)^i}\bigl(y | \omega^{2i+1} b\bigr).
\end{align}
\end{Theorem}

Note that all the factors \smash{$\tilde{g}^{(k)}_{(n-i-1)^i}\bigl(y|\omega^{2i+1}b\bigr)$}
are expressed as determinants because $(n-i-1)^i$ is $k$-small.

\subsection{Organization}
The paper is organized as follows. In Section~\ref{sec:RTL}, we first define the relativistic Toda lattice and explain
the presentation for $QK_T(\SL_n(\C)/B)$ due to Maeno, Naito, and Sagaki~\cite{MNS1,MNS2}.
We also review the centralizer family $\tilde{\mathscr{Z}}$ given in~\cite{ISY}. We introduce the $\tau$-functions as elements of $\mathscr{O}\bigl(\tilde{\mathscr{Z}}\bigr)$.
At the last part of this section we give a~detailed review of the construction of $\Phi_n$ given in~\cite{IIM}. In Section~\ref{sec:auto}, we collect results on some automorphisms used in the rest of the paper.
In Section~\ref{sec:Kk}, we review the definition and basic results on the $K$-theoretic double $k$-Schur functions.
In Section~\ref{sec:Groth}, we introduce the quantum double Grothendieck polynomials and prove Theorem~\ref{thm:main}.
As an application we prove a~formula giving a~relation between \smash{$g_\la^{(k)}(y|b)$} and \smash{$\tilde{g}_\la^{(k)}(y|b)$}.
In Section~\ref{sec:det}, we prove determinantal formulas for the $K$-theoretic double $k$-Schur functions
associated with $k$-small partitions.
In Section~\ref{sec:krect}, we prove the $k$-rectangle factorization property.
In Section~\ref{sec:max factor}, we prove Theorem~\ref{thm:max factor}.
In Appendix~\ref{sec:Hirotas_bilinear_form}, we discuss the meaning of the automorphism $\sigma$
in the context of discrete integrable systems.
In Appendix~\ref{sec:F}, we prove a~formula for the conserved quantities of the relativistic Toda lattice.
In Appendix~\ref{sec:nil-Hecke_on_QK}, we discuss the affine $K$-nil-Hecke action on the quantum $K$-theory ring and give a~proof of Proposition~\ref{prop:Quantum_D(1)}.

\subsection*{List of symbols}
\begin{itemize}[itemsep=0pt]
\item $I$, $\Iaf=I\cup\{0\}$: (affine) Dynkin index set, Section~\ref{sec:def_K-k},
 \item $Q_i$: Novikov variables, Section~\ref{sec:QK},
 \item ${\rm e}^{a_i}$, $b_i$: the equivariant parameters, Section~\ref{sec:double_K-k}, \eqref{eq:b_i},
 \item $\sigma_i$, $\tau_i$: the tau-functions,~\eqref{eq:def_of_tau}, \eqref{eq:def_of_sigma_i},
 \item $T_i$, $D_i$: the Demazure operators,~\eqref{eq:def_Demazure},
 \item $\Omega(b_i|y)$: \eqref{eq:s0_Omega},
 \item $\de_i(y)$: \eqref{eq:def de},
 \item $\co_\la(y)$: \eqref{eq:def xi},
 \item $\Phi_n$, $\tilde{\Phi}_n$: the $K$-Peterson maps,~\eqref{eq:Phi}, \eqref{eq:tPhi},
 \item $\iota$: an involution, Section~\ref{sec:iota on OZ},
 \item $T_i^Q$, $D_i^Q$: Demazure operators on the quantum $K$-ring,~\eqref{eq:DiQ,TiQ}, \eqref{eq:D0Q},
 \item $F_i(z,Q)$: the conserved quantities of the relativistic Toda lattice,~\eqref{eq:polynomialF},
 \item \smash{$F_j^{(i)}$}: \eqref{eq:Fij},
 \item \smash{$g_{x}^{(k)}(y|b)$}, \smash{$\tilde{g}_{x}^{(k)}(y|b)$}: $K$-theoretic double $k$-Schur functions,~\eqref{eq:def:tg,g},
 \item $\mathscr{Z}$, $\tilde{\mathscr{Z}}$: centralizer families, Section~\ref{ssec:tZ},
 \item $\mathscr{Z}^\circ$, $\tilde{\mathscr{Z}}^\circ$: open parts of centralizer families, Section~\ref{sec:tau},
 \item $Z=(z_{ij})$: matrix of coordinate functions of $\tilde{\mathscr{Z}}$, \eqref{eq:def-cent},
 \item $A$: \eqref{eq:def-cent},
 \item $C_A$: companion matrix of $A$, \eqref{eq:C_A},
 \item $P$: transition matrix,
 from $A$ to $C_A$, \eqref{eq:P},
 \item $L$ (and $M$, $N$): Lax matrix, Section~\ref{sec:Rel Toda},
 \item $\mathscr{Y}$, $\mathscr{Y}^\circ$: family of isolevel sets and its open part, Section~\ref{sec:QK}, Remark~\ref{rem:QK=O(Y)},
 \item $c_i$, \smash{$c_i^{(j)}$}: Section~\ref{sec:ci},
 \item \smash{$\mathfrak{G}_w^Q(z|\eta)$}: quantum double Grothendieck polynomials, Section~\ref{sec:G^Q(z|eta)},
 \item $\fac_i$: \eqref{eq:def phi},
 \item $z_i$: variables of $\mathfrak{G}_w^Q(z| \eta)$, Section~\ref{sec:G^Q(z|eta)},
 \item $\eta_i$: equivariant parameters of $\mathfrak{G}_w^Q(z| \eta)$, \eqref{eq:eta-b},
 \item $\sigma$: automorphism, Section~\ref{sec:sigma},
 \item $\omega$: a~cyclic permutation, Section~\ref{sec:def_K-k},
 \item $\omega_k$: involution,~\eqref{eq:omega_k},
 \item $R_i$: $k$-rectangle, the partition $(i)^{n-i}$,
 \item $\nu_n$: the maximal $k$-irreducible $k$-bounded partition, Section~\ref{ssec:fac max},
 \item $M_\la$: \eqref{eq:matrix_M},
 \item \smash{$\mathscr{P}^{(k)}$}: set of $k$-bounded partitions, Section~\ref{sec:double_K-k},
 \item $W_G$: Weyl group,
 \item $\hat{W}_G$: affine Weyl group, Section~\ref{sec:def_K-k},
 \item $\hat{W}_G^0$: mimimal-length coset representatives for $\hat{W}_G/W_G$.
\end{itemize}

\section{Relativistic Toda lattice and the centralizer family}\label{sec:RTL}
The aim of this section is to explain how the map $\Phi_n$ is
obtained by solving the relativistic Toda lattice.
Taking into account
recent developments,
we will review the construction in~\cite{IIM}.

\subsection{Relativistic Toda lattice}\label{sec:Rel Toda}
The \emph{relativistic Toda lattice equation}
was introduced by Ruijsenaars~\cite{Rui}.
In this paper, we start from the Lax equation due to Suris~\cite{Suris}.
Let
\begin{align}\label{eq:def_of_M_and_Q}
M=\begin{pmatrix}
z_1 & -1 & 0 & \cdots &0\\
 & z_2 & -1 & \ddots &\vdots\\
 & & \ddots & \ddots & 0\\
 & & &\ddots & -1\\
 & & & & z_n
\end{pmatrix},\quad
N=\begin{pmatrix}
1 & 0 & 0 & \cdots &0\\
-Q_1z_1 & 1 & 0 & \ddots &\vdots\\
 & -Q_2z_2 & \ddots & \ddots & 0\\
 & & \ddots&1& 0\\
 & & &-Q_{n-1}z_{n-1} & 1
\end{pmatrix}
\end{align}
and $L:=MN^{-1}$, which we call the Lax matrix.
We consider the system of partial differential equations
\begin{equation}\label{eq:rel_Toda}
{\partial L}/{\partial t_i}=\big[L,\bigl(L^i\bigr)_{<}\big]\qquad\text{for}\quad 1\le i\le n-1,
\end{equation}
where $\bigl(L^i\bigr)_{<}$ is the
strictly lower triangular part of $L^i$.
The equation is a~group version of the finite open Toda lattice.

\subsection[Quantum K-ring of the flag variety]{Quantum $\boldsymbol{K}$-ring of the flag variety}\label{sec:QK}
The integrable system is related to the quantum
$K$-ring of the flag variety $\SL_n(\C)/B$, analogously to a~result of
Givental and Kim~\cite{GK} for the quantum cohomology ring.

The conserved quantities $F_i(z,Q)$
of the system~\eqref{eq:rel_Toda} are given by
\[
\det(\zeta E-L)=\sum_{i=0}^n(-1)^i F_i(z,Q)\zeta^{n-i}.
\]
Explicitly, we have (see Appendix~\ref{sec:F})
\begin{equation}\label{eq:polynomialF}
F_i(z,Q)=
\sum_{\substack{J\subset \{1,\dots,n\}\\ |J|=i}}
\prod_{{j\in J,\, j+1\notin J}}(1-Q_j)\prod_{j\in J}z_j.
\end{equation}

Let $R(T)[\![Q]\!]$ denote the
ring of formal power series in the variables $Q_1,\dots,Q_{n-1}$ with coefficients in $R(T)$.
The quantum $K$-ring $QK_T(\SL_n(\C)/B)$ is a~commutative $R(T)[\![Q]\!]$-algebra.
\begin{Theorem}[Maeno, Naito, Sagaki~\cite{MNS1}]\label{thm:MNS1}
There exists an isomorphism of $R(T)[\![Q]\!]$-algebras
\[
QK_T(\SL_n(\C)/B)\cong
\Rep[\![Q]\!][z_1,\dots,z_n]/\mathscr{I}_n^Q,
\]
where the ideal $\mathscr{I}_n^Q$ is generated by the elements
\begin{equation}\label{eq:conserved}
F_i(z,Q)-e_i({\rm e}^{-a_1},\dots,{\rm e}^{-a_n})\qquad \text{for}\quad 1\le i\le n. \end{equation}
\end{Theorem}
\begin{Remark}
Our $a_i$ is $-\epsilon_i$ in~\cite{MNS1}.
\end{Remark}In this geometric
context, $z_i$ is related to the class of
the universal line bundle (see~\cite{MNS1} for details),
and $Q_i$ $(1\le i\le n-1)$ are the Novikov variables.

\begin{Remark}\label{rem:KPSZ}
It follows from \cite[Corollary~2]{GL} that the elements in~\eqref{eq:conserved} are contained in the defining ideal of $QK_T(\SL_n(\C)/B)$ (see also~\cite{ACT, MNS1} for details).
We also note that Koroteev--Puskar--Smirnov--Zeitlin~\cite{KPSZ} proved that the elements in~\eqref{eq:conserved} generate the defining ideal of $
QK_T^{\mathrm{QM}}\bigl(\SL_n(\C)/B\bigr)$,
which is the limit of quasimap quantum $K$-ring of $T^*(\SL_n(\C)/B)$ at $\hbar = \infty$.
Moreover, Huq--Kuruvilla~\cite{HK} established that this limit is isomorphic to $QK_T(\SL_n(\C)/B)$ (more generally, for partial flag varieties).
\end{Remark}

Let $\mathscr{Y}$ be the affine subscheme
of $\mathbb{A}^{2n-1}$ with coordinates
$Q_1,\dots,Q_{n-1}$, $z_1,\dots,z_n$
whose defining ideal is generated by the polynomials~\eqref{eq:conserved}.
$\mathscr{Y}$ is a~scheme over $T$ whose fibers are the isolevel sets.
The coordinate ring, i.e., the ring of regular functions $\mathscr{O}(\mathscr{Y})$ is
considered to be a~polynomial version \smash{$QK_T^{\mathrm{pol}}(\SL_n(\C)/B)$} of $QK_T(\SL_n(\C)/B)$.

\begin{Remark}
\label{rem:QK=O(Y)}
Let us denote the ideal of $\Rep[z,Q]$
generated by the polynomials~\eqref{eq:conserved}
by $\mathscr{I}_n^{Q,\mathrm{pol}}$.
The quotient ring
\smash{$\mathscr{O}(\mathscr{Y})=\Rep[z,Q]/\mathscr{I}_n^{Q,\mathrm{pol}}$} is
\smash{$QK_T^{\mathrm{pol}}(\SL_n(\C)/B)$} in the notation used in the introduction.
Although there are subtleties in relating \smash{$QK_T^{\mathrm{pol}}(\SL_n(\C)/B)$}
to $QK_T(\SL_n(\C)/B)$,
\smash{$QK_T^{\mathrm{pol}}(\SL_n(\C)/B)_Q=\mathscr{O}(\mathscr{Y}^\circ)$}
is actually a~subring of $QK_T(\SL_n(\C)/B)_Q$. Furthermore, the residue class of \smash{$\mathfrak{G}_w^Q$} in \smash{$QK_T^{\mathrm{pol}}(\SL_n(\C)/B)_Q$} corresponds to
$\mathscr{O}^w$ in $QK_{T}(\SL_n(\C)/B)$. For these facts, see \cite[Sections~4.1 and 4.2]{IIN}; it is straightforward to have the equivariant version based on \cite[Remark~6.3]{MNS1}.
\end{Remark}

If we consider
the special fiber $\mathscr{Y}_{\mathrm{uni}}$
corresponding to the case when the Lax matrix is \emph{unipotent}, and the corresponding
centralizer $\mathscr{Z}_{\mathrm{uni}}$, we obtain the non-equivariant $K$-Peterson isomorphism
$QK(\SL_n(\C)/B)_\loc^{\mathrm{pol}}\cong K_*(\Gr_{\SL_n})_\loc$ studied in~\cite{IIM}.
In fact, the element of $\mathscr{Z}_{\mathrm{uni}}$ is of the form
\[
\begin{pmatrix}
1 & h_1 & h_2 & \cdots &h_{n-1}\\
 & 1 & h_1 & \ddots &\vdots\\
 & & \ddots & \ddots & h_2\\
 & & &\ddots & h_1\\
 & & & & 1
\end{pmatrix},
\]
so $\mathscr{O}(\mathscr{Z}_\mathrm{uni})$ is a~polynomial ring of $(n-1)$ variables, which can be
identified with
the non-equivariant $K$-homology ring
$K_*(\Gr_{\SL_n})$ studied by
Lam, Schilling, and Shimozono.
If we identify the generators
of the polynomial ring with the complete symmetric functions
$h_1,\dots,h_{n-1}$, the Schubert structure sheaves are
given by the $K$-theoretic closed $k$-Schur functions. See~\cite{IIN} for a~more precise correspondence.

\subsection[Centralizer family Z]{Centralizer family $\boldsymbol{\tilde{\mathscr{Z}}}$}
\label{ssec:tZ}
We construct the solutions of~\eqref{eq:rel_Toda} from a~centralizer family $\mathscr{Z}$ defined below.
Consider the matrix equation
\begin{equation}
\label{eq:def-cent}
[A,Z]=0, \qquad
A=\begin{pmatrix}
{\rm e}^{-a_1} & -1 & 0 & \cdots &0\\
 & {\rm e}^{-a_2} & -1 & \ddots &\vdots\\
 & & \ddots & \ddots & 0\\
 & & &\ddots & -1\\
 & & & & {\rm e}^{-a_n}
\end{pmatrix},
\end{equation}
where $Z$ is an upper triangular
matrix with the entries $z_{ij}$ for $1\le i\le j\le n$. We assume $z_{11}\cdots z_{nn}\ne 0$.
$[A,Z]=0$ is equivalent to the equations
\begin{equation}
(b_i-b_j){z}_{ij}={{z}_{i,j-1}-{z}_{i+1,j}}\qquad \text{for} \quad1\le i\le j\le n. \label{eq:cent}
\end{equation}

Let us denote the affine variety over $T$ defined by these equations by $\tilde{\mathscr{Z}}$. So we have
\begin{gather*}
\mathscr{O}\bigl(\tilde{\mathscr{Z}}\bigr)=
\Rep\big[z_{ii}^{\pm 1} \ (1\le i\le n),\,
z_{ij} \ (1\le i< j\le n)
\big]/I,\\
I=\langle
(b_i-b_j)z_{ij}-z_{i,j-1}+z_{i+1,j}
\mid 1\le i< j\le n
\rangle.
\end{gather*}
There is a~natural
$\C^\times$-action
on $\tilde{\mathscr{Z}}$ by scalar
multiplication $z_{ij}\mapsto
cz_{ij}$ $(c\in \C^\times)$.
The variety $\mathscr{Z}$ is defined to be
$\tilde{\mathscr{Z}}/\C^\times$.
Thus the coordinate ring
$\mathscr{O}({\mathscr{Z}})$
is the $\Rep$-subalgebra of $\mathscr{O}\bigl(\tilde{\mathscr{Z}}\bigr)$ generated by
$z_{ij}/z_{11}$
$(1\le i\le j\le n)$.
$\mathscr{Z}$ is a~closed subscheme of
$T\times B^\vee$, where
$B^\vee$ is the
Borel subgroup of $\PGL_n(\C)$.

\begin{Remark}
In~\cite{ISY},
$\mathscr{Z}$ and
$\tilde{\mathscr{Z}}$
are denoted by $\mathscr{Z}_{\PGL_n}$
and $\mathscr{Z}_{\GL_n}$
respectively.
\end{Remark}

Let $\Rep^\reg$ be the localization of $\Rep$ by
the multiplicative set generated by $1-{\rm e}^{a_i-a_j}$ $(i\ne j)$.
Set \smash{$\mathscr{O}\bigl(\tilde{\mathscr{Z}}\bigr)^\reg
:=\Rep^\reg\otimes_{\Rep}\mathscr{O}\bigl(\tilde{\mathscr{Z}}\bigr)$}.
We often use the following.
\begin{Proposition}\label{prop:OZ_Delta}
$\mathscr{O}\bigl(\tilde{\mathscr{Z}}\bigr)^\reg$
is generated by $z_{ii}$ $(1\le i\le n)$ as
an~$\Rep^\reg$-algebra.
\end{Proposition}

\subsection[tau-functions]{$\boldsymbol{\tau}$-functions}\label{sec:tau}
We will explain that a~Zariski open subset $\mathscr{Z}^\circ$ of $\mathscr{Z}$ is isomorphic to
$\mathscr{Y}^\circ$, the compliment of the divisors defined by $Q_i$ $(1\le i\le n-1)$, whereas
$\mathscr{Z}^\circ$ is the complement of the divisor defined by the so-called $\tau$-\emph{functions}.

Let \smash{${\rm e}^{(m)}_i:=e_i({\rm e}^{-a_1},{\rm e}^{-a_2},\dots,{\rm e}^{-a_m})$} denote the $i$-th elementary symmetric polynomial in $m$ variables ${\rm e}^{-a_1},{\rm e}^{-a_2},\dots,{\rm e}^{-a_m}$.
The companion matrix $C_A$ of $A$ is described as
\begin{equation}
\label{eq:C_A}
C_A:=P^{-1}AP
=\begin{pmatrix}
0& -1 \\
&0& -1 \\
& & \ddots&\ddots \\
&& & 0&-1\\
{\rm e}^{(n)}_n& {\rm e}^{(n)}_{n-1}&\cdots &{\rm e}^{(n)}_2&{\rm e}^{(n)}_1\\
\end{pmatrix},
\end{equation}
where
\begin{equation}
\label{eq:P}
P=
\begin{pmatrix}
1\\
{\rm e}^{(1)}_1 & 1\\[1mm]
{\rm e}^{(2)}_2 & {\rm e}^{(2)}_1 & 1\\[1mm]
{\rm e}^{(3)}_3 & {\rm e}^{(3)}_2 & {\rm e}^{(3)}_1&1\\[1mm]
\vdots & \vdots & \vdots & & \ddots\\
{\rm e}^{(n-1)}_{n-1}& {\rm e}^{(n-1)}_{n-2}&\cdots &\cdots&{\rm e}^{(n-1)}_1 & 1
\end{pmatrix}.
\end{equation}
For later use, we record the following formula:
\[
 \bigl(P^{-1}\bigr)_{ij}=(-1)^{i-j}h_{i-j}({\rm e}^{-a_1},\dots,{\rm e}^{-a_{j}}).
\]

Let $I$, $J$ be subsets of $[1,n]:=\{1,\dots,n\}$,
with $I=\{i_1,\dots,i_p\}$ and $J=\{j_1,\dots,j_q\}$.
Let
$X^{J}_{I}$ denote the submatrix of $X$ consisting of its $i_1,\dots,i_p$-th rows and $j_1,\dots,j_q$-th columns.
\begin{Definition}
For $1\le i\le n$, define $\tau_i, \sigma_i\in \O\bigl(\tilde{\mathscr{Z}}\bigr)$ by
\begin{gather}\label{eq:def_of_tau}
\tau_i=
\det(ZAP)^{[1,i]}_{[1,i]},\\
\label{eq:def_of_sigma_i}
\sigma_i
=
\det(ZP)^{[1,i]}_{[1,i]},
\end{gather}
and $\tau_0=\sigma_0=1$.
\end{Definition}
In particular, we have
$
\sigma_n=\det(Z)=z_{11}\cdots z_{nn}$, $
\tau_n={\rm e}^{a_1+\dots+a_n}\sigma_n$.
The open set $\tilde{\mathscr{Z}}^\circ$
of $\tilde{\mathscr{Z}}$ is
the complement of the closed set defined by
$\tau_i= 0$, $\sigma_i=0$ ($1\le i\le n-1$).

We often use the following results on
the determinants.
\begin{Proposition}\label{prop:Noumi}
Let $A$, $B$ be square matrices of size $n$, and $1\le i\le n$.
If $A$ is lower triangular or $B$ is upper triangular, then
\[
\det(AB)^{[1,i]}_{[1,i]}=\det A^{[1,i]}_{[1,i]}
\det B^{[1,i]}_{[1,i]}.
\]
\end{Proposition}

\begin{Proposition}\label{prop:det AB partial}
Let $A$, $B$ be square matrices of size $n$ such that $B$ is invertible, and $1\le i\le n-1$.
\[
\det(AB)^{[i+1,n]}_{[i+1,n]}
=\det\left(
\begin{matrix}
\bigl(B^{-1}\bigr)^{[1,n]}_{[1,i]}\\\hline
A^{[1,n]}_{[i+1,n]}
\end{matrix}
\right)\cdot \det B.
\]
\end{Proposition}

\subsection[Construction of Phi\_n]{Construction of $\boldsymbol{\Phi_n}$}
\begin{Theorem}[{\cite[Section~3]{IIM}}]By the map defined by~\eqref{eq:Phi}, we have
an isomorphism
$\mathscr{Y}^\circ\cong \mathscr{Z}^\circ$
of varieties over $T$.
\end{Theorem}

An element of $\tilde{\mathscr{Z}}$
is an algebraic family $\{Z_t\}$ of invertible
upper triangular matrices
parame\-trized by $t\in T$ such that $Z_t$ commute with $A_t$.
We denote the family $\{Z_t\}$ simply by $Z$.
If we assume~${Z\in \tilde{\mathscr{Z}}^\circ}$,
there exists
an upper triangular matrix $R$
and a~unipotent lower triangular matrix~$U$ (see \cite[Proposition 3.2]{IIM}), both defined as families over $T$, such that \begin{equation}\label{eq:Gauss_decomp}
P^{-1}ZAP=U^{-1}R.
\end{equation}
Because the left-hand side of~\eqref{eq:Gauss_decomp} commutes with $C_A$, we have
\begin{equation}\label{eq:AZ_to_L} UC_AU^{-1}=RC_AR^{-1}.
\end{equation}
Let $L$ denote the matrix~\eqref{eq:AZ_to_L}.
If we replace $Z$ by $cZ$ with $c\in \C^\times$,
then $R$ becomes $cR$ and we obtain the same matrix $L$.
We apply the Gauss decomposition to
$L$ as $L=MN^{-1}$. Then $M$, $N$ are matrices of the forms given in~\eqref{eq:def_of_M_and_Q} for unique $z_1,\dots,z_n, Q_1,\dots,Q_{n-1}$ (see~\cite[Section~3.2]{IIM}).
Thus we obtain
functions
$z_i$, $Q_i$
on $\mathscr{Z}^\circ=\tilde{\mathscr{Z}}^\circ/\C^{\times}$.
In view of~\eqref{eq:AZ_to_L},
we have~${\det(\zeta E-L)
=\det(\zeta E-C_A)}$, which means that
the regular functions $z_1,\dots,z_n$, $Q_1,\dots,\allowbreak Q_{n-1}$
satisfy the definiting equation for $\mathscr{Y}$.
We can also check $Q_i\ne 0$
for $1\le i\le n-1$ (see \cite[Section~3.5]{IIM}).
In this way, we have a~ring homomorphism
$\Phi_n\colon \mathscr{O}(\mathscr{Z}^\circ)\rightarrow\mathscr{O}(\mathscr{Y}^\circ)$, which is naturally a~homomorphism of $\Rep$-algebras.

The explicit form of $\Phi_n$
is determined as follows.
By abuse of notation
we simply denote $\Phi_n(z_i)$, $\Phi_n(Q_i)$ by $z_i$, $Q_i$.
Let $r_{ij}$ be the $(i,j)$-th entry of $R$.
Then, from~\eqref{eq:Gauss_decomp} and the Cauchy--Binet formula, we obtain
\begin{equation}\label{eq:r_to_minor}
\begin{aligned}
r_{11}r_{22}\cdots r_{ii}
&=
\det R^{[1,i]}_{[1,i]}=
\det\bigl(UP^{-1}ZAP\bigr)^{[1,i]}_{[1,i]}
=\det(ZAP)^{[1,i]}_{[1,i]}=\tau_i,
\end{aligned}
\end{equation}
which implies $r_{ii}=\tau_{i}/\tau_{i-1}$.
Comparing the $(i+1,i)$-th entries on both sides of
$
NM^{-1}=RC_A^{-1}R^{-1}
$, which is derived from~\eqref{eq:AZ_to_L}, we obtain
\begin{equation}\label{eq:Q_to_tau}
Q_i=\frac{r_{i+1,i+1}}{r_{ii}}=\frac{\tau_{i+1}\tau_{i-1}}{\tau_i^2}.
\end{equation}

On the other hand, from~\eqref{eq:def_of_M_and_Q} and~\eqref{eq:AZ_to_L}, $z_i$ is expressed as
\[
z_{i}
=\frac{\det\bigl(L^{-1}\bigr)^{[1,i-1]}_{[1,i-1]}}{\det\bigl(L^{-1}\bigr)^{[1,i]}_{[1,i]}}
=\frac{\det\bigl(UC_A^{-1}U^{-1}\bigr)^{[1,i-1]}_{[1,i-1]}}{\det\bigl(UC_A^{-1}U^{-1}\bigr)^{[1,i]}_{[1,i]}}
=
\frac{\det\bigl(C_A^{-1}U^{-1}\bigr)^{[1,i-1]}_{[1,i-1]}}{\det\bigl(C_A^{-1}U^{-1}\bigr)^{[1,i]}_{[1,i]}}.
\]
As $C_A^{-1}U^{-1}
=P^{-1}ZPR^{-1}$, we have
\begin{align*}
\det\bigl(C_A^{-1}U^{-1}\bigr)^{[1,i]}_{[1,i]}
&
=\det\bigl(P^{-1}ZPR^{-1}\bigr)^{[1,i]}_{[1,i]}
=\frac{\det\bigl(P^{-1}ZP\bigr)^{[1,i]}_{[1,i]}}{r_{11}\cdots r_{ii}}
=\frac{\det(ZP)^{[1,i]}_{[1,i]}}{r_{11}\cdots r_{ii}}
=\frac{\sigma_i}{\tau_i},
\end{align*}
which gives
\begin{equation}\label{eq:z_to_tau}
z_i
=\frac{\tau_i\sigma_{i-1}}{\tau_{i-1}\sigma_i}.
\end{equation}

We can construct the inverse of
$\Phi_n$ (see \cite[Section~3.4]{IIM}).
Here we provide an expression for~$\Phi_n^{-1}$.
From~\eqref{eq:Gauss_decomp}, we have
$
P^{-1}ZP=U^{-1}R C_A^{-1}$.
The entries of $U^{-1}$ can be expressed as
polynomials in $z_i$, $Q_i$, while the entries of $R$ can be written as Laurent polynomials in $z_i$, $Q_i$.
This gives rise to an expression
for the entries of $Z$ as Laurent polynomials in $z_i$, $Q_i$, which are considered as elements in $\mathscr{O}(\mathscr{Y}^\circ)$.
\begin{Example} For $n=2$, we have
\[
Z=
\frac{1}{z_1z_2Q_1}
\begin{pmatrix}
z_2-{\rm e}^{-a_1}&1\\
0&z_2-{\rm e}^{-a_2}
\end{pmatrix}.
\]
For $n=3$, we have
\begin{align*}
&z_{11}=(z_1z_2z_3Q_1Q_2)^{-1}\bigl(z_2z_3-{\rm e}^{-a_1}(z_2(1-Q_2)+z_3)+{\rm e}^{-2a_1}\bigr),\\
&z_{12}=(z_1z_2z_3Q_1Q_2)^{-1}\bigl(z_2(1-Q_2)+z_3-{\rm e}^{-a_1}-{\rm e}^{-a_2}\bigr),\\
&z_{13}=(z_1z_2z_3Q_1Q_2)^{-1}.
\end{align*}
The other entries are determined by
$z_{i+1,j+1}=\omega(z_{ij})$.
\end{Example}

\begin{Remark}
As a~complex manifold, $\mathscr{Z}^\circ$ is parameterized as
\begin{equation}\label{eq:t_to_Z}
Z=\exp\bigl(
At_1+A^2t_2+\cdots+A^{n-1}t_{n-1}
\bigr)\in \mathscr{Z}^\circ
\end{equation}
by the complex parameters $t_1,t_2,\dots,t_{n-1}$.
From~\eqref{eq:t_to_Z}, we deduce the differential equation of motion $\frac{\partial }{\partial t_i}Z=A^iZ$.
Then, computing the differential $\frac{\partial}{\partial t_i}\bigl(P^{-1}ZAP\bigr)$, we obtain
\[
\frac{\partial}{\partial t_i}\bigl(P^{-1}ZAP\bigr)=P^{-1}\left(\frac{\partial}{\partial t_i}Z\right)AP=P^{-1}A^iZAP=
C_A^iU^{-1}R,
\]
where, for the last equality, we used~\eqref{eq:Gauss_decomp}.
On the other hand, by using~\eqref{eq:Gauss_decomp} again, we also obtain
\[
\frac{\partial}{\partial t_i}\bigl(P^{-1}ZAP\bigr)
=
\frac{\partial}{\partial t_i}\bigl(U^{-1}R\bigr)
=U^{-1}\left(\frac{\partial}{\partial t_i}R\right)-U^{-1}\left(\frac{\partial}{\partial t_i}U\right)U^{-1}R.
\]
Comparing these equations, we have
\[
L^i=UC_A^iU^{-1}=
\left(\frac{\partial}{\partial t_i}R\right)R^{-1}-\left(\frac{\partial}{\partial t_i}U\right)U^{-1}.
\]
Since $\bigl(\frac{\partial}{\partial t_i}R\bigr)R^{-1}$ is upper triangular, we obtain \smash{$\bigl(L^i\bigr)_{<}=-\bigl(\frac{\partial}{\partial t_i}U_{t_i}\bigr)U^{-1}$} by~\eqref{eq:AZ_to_L}. Then, we have the relativistic Toda lattice~\eqref{eq:rel_Toda} as follows:
\begin{align*}
\frac{\partial}{\partial t_i}L
&
=
\frac{\partial}{\partial t_i}\bigl(UC_AU^{-1}\bigr)=
\left(\frac{\partial}{\partial t_i}U\right)C_AU^{-1}-UC_AU^{-1}
\left(\frac{\partial}{\partial t_i}U\right)U^{-1}
\\
&=
\left(\frac{\partial}{\partial t_i}U\right)U^{-1}L-L
\left(\frac{\partial}{\partial t_i}U\right)U^{-1}=\big[L,\bigl(L^i\bigr)_{<}\big].
\end{align*}
\end{Remark}

\section[K-theoretic double k-Schur functions]{$\boldsymbol{K}$-theoretic double $\boldsymbol{k}$-Schur functions}\label{sec:Kk}
In order to study
the map $\Phi_n$
in more detail,
we use a~recent
work~\cite{ISY} on a~symmetric function realization
of $K_*^T(\Gr_{\SL_n})$.
\subsection[Definition of K-theoretic double k-Schur functions]{Definition of $\boldsymbol{K}$-theoretic double $\boldsymbol{k}$-Schur functions}
\label{sec:def_K-k}
The simple roots of $G=SL_n(\C)$ are given by
$\alpha_i=a_i-a_{i+1}$ for $i\in I$, where
$I=\{1,\dots,n-1\}$.
Let $\theta=a_1-a_n$ denote the highest root, and let $\theta^\vee$ be its corresponding coroot.
We define an action of $\tilde{W}_G$ on $R(T)$, referred to as the \emph{level zero} action. The coroot lattice $Q^\lor$ acts by the identity, while
$W_G$ acts naturally. In particular, $s_0$ acts as $s_\theta$ where $\theta$ is the highest root.

We will work in the \emph{level zero} affine setting so we set $\alpha_0:=-\theta$.
Let $\tilde{W}_G$ denote the corresponding affine Weyl group, with the standard generators $s_0,s_1,\dots,s_{n-1}$.
Let $\Rep^\reg$ denote the localization of $\Rep$ by
the multiplicative set generated by $1-{\rm e}^{\alpha}$ where
$\alpha$ are any roots.
The \emph{twisted group algebra} $\Rep^\reg\big[\tilde{W}_G\big]$
is $\bigoplus_{w\in\tilde{W}_G }\Rep^\reg w$ with product defined by
\[
(f_1w_1)(f_2w_2)=(f_1w_1(f_2))(w_1w_2)
\qquad \text{for}\quad f_1,f_2\in\Rep^\reg,\quad w_1,w_2\in \tilde{W}_G.
\]
with the level zero action of
$\tilde{W}_G$ on $\Rep^\reg$.

For $i\in \tilde{I}=I\cup\{0\}$, define the \emph{ Demazure operators}
\begin{equation}
\label{eq:def_Demazure}
T_i=(1-{\rm e}^{\alpha_i})^{-1}(s_i-1),\qquad
D_i=T_i+1,
\end{equation}
which are considered as
elements of \smash{$\Rep^\reg\big[\tilde{W}_G\big]$}.

$T_i$ satisfies $T_i^2 = -T_i$ and the braid relations of type \smash{$A_{n-1}^{(1)}$}. Similarly, $D_i$ satisfies $D_i^2 = D_i$ and the braid relation of type \smash{$A_{n-1}^{(1)}$}. Then, for any $ w \in \hat{W}_G$ written as $w = s_{i_1} \cdots s_{i_l}$, the product
$
T_w = T_{i_1} \cdots T_{i_l}
$ and the product
$
D_w = D_{i_1} \cdots D_{i_l}
$ depend only on $w$ (see \cite{LSS:K}, see also~\cite{ISY}).

The $K$-theoretic nil-Hecke algebra $\mathbb{K}_G$ is defined to be the left $\Rep$-module
generated by~$D_w$ for $w\in \hat{W}_G$
(see \cite{LSS:K}, see also~\cite{ISY}).

We use notation
\begin{equation}
\label{eq:b_i}
b_i:=1-{\rm e}^{-a_i},\qquad 1\le i\le n.
\end{equation}

We consider the
ring \smash{$\hat{\Lambda}^{\Rep}$} of symmetric
formal power series
in the infinitely many variables~${y=(y_1,y_2,\dots)}$ with coefficients in
$\Rep$.
We can define an action of Demazure operators~$D_i$ $(0\le i\le n-1)$
on $\hat{\Lambda}^{\Rep}$ (see~\cite{ISY}).
For $1\le i\le n-1$, we
let $s_i$ act by exchanging
$a_i$ and~$a_{i+1}$
in the coefficients,
and $s_0$ by
the formula
\begin{equation}
\label{eq:s0_Omega}
s_0(f)=\frac{\Omega(b_1|y)}{\Omega(b_n|y)}s_\theta(f),\qquad
\Omega(b_i|y):=\frac{1}{\prod_{j=1}^{\infty}(1-b_i y_j)},
\end{equation}
where $s_\theta $ is the reflection with respect to $\theta$,
exchanging $a_1$ and $a_n$.

The finite Weyl group $W=\langle s_1,\dots,s_{n-1}\rangle$ is the symmetric group $S_n$.
Let \smash{$\hat{W}_G^0$} denote the set of minimal-length coset representatives of \smash{$\hat{W}_G/W_G$}. The set \smash{$\hat{W}_G^0$} naturally indexes the set of the Schubert classes of \smash{$K_*^T(\Gr_{G})$},
and we refer to an element of \smash{$\hat{W}_G^0$} as an affine Grassmann element.
For each
\smash{$x\in \hat{W}_G^0$}
\cite{ISY},
the \emph{double $K$-theoretic $k$-Schur functions}
are defined
as
\begin{equation}
\label{eq:def:tg,g}
\tilde{g}_x^{(k)}(y|b):=D_x (1),\qquad
g_x^{(k)}(y|b):=T_x (1).
\end{equation}
The $\Rep$-span
\[
\hat{\Lambda}^{\Rep}_{(n)}:=\bigoplus_{x\in \hat{W}_G^0}
\Rep \tilde{g}_x^{(k)}(y|b)
\]
is an $\Rep$-subalgebra of
\smash{$\hat{\Lambda}
_\Rep$} which is
isomorphic to
\smash{$K_*^T(\Gr_{\SL_n})$} and the Schubert structure~sheaf $\mathscr{O}_x$
is identified with \smash{$\tilde{g}_x^{(k)}(x|b)$}.
Let $\omega\in S_n$ be the cyclic permutation sending $i$ to $i+1$ with $n+1=1$ by convention.

\begin{Theorem}[\cite{ISY}]
\label{thm:ISY}
There are
isomorphisms $\beta$, $\kappa$ of $R(T)$-algebras
\[
\mathscr{O}(\mathscr{Z})
\overset{\beta}{\longrightarrow }
 \hat{\Lambda}_{(n)}^\Rep
 \overset{\kappa}{\longrightarrow}
 K_*^T(\Gr_{\SL_n}),
\]
such that
\begin{gather}
\beta(z_{ij}/z_{11})= {\rm e}^{a_i+\cdots+a_{j-1}}\Omega(b_i|y)\Omega(b_1|y)^{-1}g_{\rho_{j-i}}\big(y|\omega^i(b)\big),\nonumber
\\
\kappa(\tilde{g}_{x}^{(k)}(y|b))
= \mathscr{O}_x,
\qquad
\kappa(g_{x}^{(k)}(y|b))=\mathscr{I}_{x},\label{eq:zij_in_KkSchur}
\end{gather}
where $\rho_l=s_{l-1}\cdots s_{1}s_{0}$ for $1\le l\le n-1 $.
\end{Theorem}

\begin{Remark}
In~\cite{ISY}, it is proven that
each of the three $R(T)$-algebras has a~natural $\Rep$-Hopf-algebra structure.
\end{Remark}

\begin{Remark}
It is natural to define a~map
\smash{$\beta\colon \mathscr{O}\bigl(\tilde{\mathscr{Z}}\bigr)
\rightarrow \hat{\Lambda}^\Rep$}. We have
\begin{equation}
\beta(z_{ij})= {\rm e}^{a_i+\cdots+a_{j-1}}\Omega(b_i|y)g_{\rho_{j-i}}\big(y|\omega^i(b)\big),
\label{eq:zij_in_KkSchur_GL}
\end{equation}
where $T$ is the maximal torus of $\GL_n(\C)$. See \cite[Theorem 1.2]{ISY} for more details.
\end{Remark}

\subsection{Preliminaries on Demazure operators}
We collect some properties of $T_i$ and $D_i$ which are used in Section~\ref{sec:det_formula}.
Let $\iota$ be the involution on \smash{$\hat{\Lambda}^{R(T)}$} such that
\[
\iota(h_i(y))=\sum_{r=0}^{i-1}\binom{i-1}{r}e_{r+1}(y),\qquad
\iota({\rm e}^{a_i})={\rm e}^{-a_{n-i+1}}.
\]
\begin{Proposition}\label{prop:iota Dem}
For $i\in I$, $\iota\circ T_i=T_{n-i}\circ \iota$, $
\iota\circ D_i=D_{n-i}\circ \iota$.
\end{Proposition}
\begin{proof}
Straightforward.
\end{proof}

Let $\omega_k$ be the group automorphism of $\hat{W}_G$ such that
\begin{equation}
\label{eq:omega_k}
\omega_k(s_i)=s_{-i},
\end{equation}
with indices taken mod $n\Z$.
\begin{Proposition}[{\cite[Proposition 3.14]{ISY}}]
\label{prop:iota $K$-k}
We have
$
\iota\bigl(\tilde{g}_x^{(k)}\bigl(y|\omega^i b\bigr)\bigr)=
\tilde{g}^{(k)}_{x^{\omega_k}}\bigl(y|\omega^{-i}b\bigr)$.
\end{Proposition}

We define elements $T_{\theta}$ and $D_\theta$ of the twisted group algebra $\Rep^\reg\big[\tilde{W}_G\big]$:
\begin{equation}
\label{eq:Ttheta}
T_\theta:=\omega^{-i}\circ T_i \circ \omega^i=\frac{s_\theta-1}{1-{\rm e}^{-\theta}},\qquad
D_\theta:=\omega^{-i}\circ D_i \circ \omega^i=T_\theta+1.
\end{equation}

\begin{Remark}
For any $\alpha\in Q$, an element $T_\alpha$ of $\Rep^\reg\big[\tilde{W}_G\big]$ is defined. See \cite[Section~2.2.4]{ISY} for more details.
\end{Remark}

\begin{Lemma}\label{lem:T_and_D}
For $i\in I$,
$T_{i}\circ {\rm e}^{-a_{i+1}}={\rm e}^{-a_{i}}\circ D_{i}$, and
$T_{\theta} \circ {\rm e}^{-a_1}={\rm e}^{-a_n}\circ D_{\theta}$.
\end{Lemma}
\begin{proof}
The lemma is shown by direct calculations:
$T_{j}({\rm e}^{-a_{j+1}}f)=T_{j}({\rm e}^{-a_{j+1}})f+s_{j}({\rm e}^{-a_{j+1}})T_{j}(f)\allowbreak
={\rm e}^{-a_{j}}f+{\rm e}^{-a_{j}}T_{j}(f)={\rm e}^{-a_{j}}D_{j}(f)$.
\end{proof}

\begin{Proposition}\label{prop:T_0}
 We have
 \[
 T_0=\Omega(b_1|y)\circ T_{\theta}\circ \Omega(b_1|y)^{-1}
 =\Omega(b_n|y)^{-1}\circ T_{\theta}\circ \Omega(b_n|y).
 \]
\end{Proposition}
\begin{proof}
The first equality is shown as follows:
\begin{gather*}
\Omega(b_1|y)\circ T_{\theta}\circ \Omega(b_1|y)^{-1}(f)\\
\qquad
=
\Omega(b_1|y)
T_{\theta}\left(
\frac{f}{\Omega(b_1|y)}
\right)
=
\frac{1-b_n}{b_n-b_1}
\Omega(b_1|y)
\left(
\frac{f}{\Omega(b_1|y)}
-
\frac{s_\theta(f)}{\Omega(b_n|y)}
\right)\\
\qquad
=
\frac{1-b_n}{b_n-b_1}
\left(
f
-
\frac{\Omega(b_1|y)}{\Omega(b_n|y)}s_{\theta}(f)
\right)=T_0(f).
\end{gather*}
The second equality follows from the fact that $\Omega(b_1|y)\Omega(b_n|y)$ commutes with $T_{\theta}$.
\end{proof}

\begin{Lemma}\label{lem:T_on_h}
For
$1\le i\le n-1$ and $j-i\le n-1$, and $m\ge 0$,
\begin{equation}
D_{i}
h_m({\rm e}^{-a_{i+1}},{\rm e}^{-a_{i+2}},\dots,{\rm e}^{-a_j})
=
h_m({\rm e}^{-a_{i}},{\rm e}^{-a_{i+1}},{\rm e}^{-a_{i+2}},\dots,{\rm e}^{-a_j}).
\label{eq: Di hm}
\end{equation}
For $1\le j\le n-1$, and $m\ge 0$,
\begin{equation}
D_{\theta}
h_m({\rm e}^{-a_{1}},{\rm e}^{-a_{2}},\dots,{\rm e}^{-a_j})
=
h_m({\rm e}^{-a_{n}},{\rm e}^{-a_{1}},{\rm e}^{-a_{2}},\dots,{\rm e}^{-a_j}).
\label{eq: D theta hm}
\end{equation}
\end{Lemma}
\begin{proof}
Note
\begin{equation*}
D_i (1-{\rm e}^{-a_{i+1}}u)^{-1}=
(1-{\rm e}^{-a_{i}}u)^{-1}
(1-{\rm e}^{-a_{i+1}}u)^{-1}.
\end{equation*}
As $\prod_{s=i+2}^j(1-{\rm e}^{-a_s}u)^{-1}$ is invariant under $s_i$, we obtain
$D_i\prod_{s=i+1}^j(1-{\rm e}^{-a_s}u)^{-1}=
\prod_{s=i}^j(1-{\rm e}^{-a_s}u)^{-1}$, which implies~\eqref{eq: Di hm}.
The proof of~\eqref{eq: D theta hm} is given similarly.
\end{proof}

\section{Some automorphisms}\label{sec:auto}
For later use, we introduce
some automorphisms of the algebras $\mathscr{O}\bigl(\tilde{\mathscr{Z}}\bigr), \hat{\Lambda}^{\Rep}$
and study their properties.
\subsection[Automorphism sigma]{Automorphism $\boldsymbol{\sigma}$}\label{sec:sigma}

We define an automorphism $\sigma$
of $\mathscr{O}\bigl(\tilde{\mathscr{Z}}\bigr)$ as an $\Rep$-algebra
by
\begin{align}\label{eq:def_of_sigma}
\sigma(Z)&=ZA^{-1}.
\end{align}
Explicitly,
for $1\le i\le j\le n$,
we have
\[
\sigma(z_{ij})={\rm e}^{a_j}z_{ij}+
{\rm e}^{a_{j-1}+a_j}z_{i,j-1}+\cdots+{\rm e}^{a_i+\cdots+a_j}
z_{ii}.
\]
In particular,
we have
\begin{equation}
\label{eq:sigma_zii}
\sigma(z_{ii})={\rm e}^{a_i}z_{ii}.
\end{equation}
From~\eqref{eq:def_of_tau} and~\eqref{eq:def_of_sigma_i}, we have
\begin{equation}\label{eq:sigma_tau_i}
\sigma(\tau_i)=\sigma_i.
\end{equation}

Via the isomorphism $\beta$
we consider the corresponding automorphism on \smash{$\hat{\Lambda}_{(n)}^\Rep$}
and denote it also by $\sigma$.
For example,
we have
\begin{equation}
\sigma(\Omega(b_i|y))={\rm e}^{a_i}\Omega(b_i|y)\label{eq:sigma(Omega)}.
\end{equation}
In the next section, we prove
\begin{equation}
\sigma(\ks_{\rho_i}^{(k)}(y|b))={\rm e}^{a_i-a_n}\kst_{\rho_i}^{(k)}(y|b).\label{eq:sigma(tg)}
 \end{equation}
Since we know
$g_{\rho_i}^{(k)}(y|0)=h_i(y)$,
we have
$
\sigma(h_i(y))=1+h_1(y)+\cdots+h_i(y)$.
In fact, $\sigma$
on~\smash{$\hat{\Lambda}^\Rep$} can be defined
as the $\Rep$-linear map
sending $f(y_1,y_2,\dots)\in \hat{\Lambda}$
to $f(1,y_1,\allowbreak y_2,\dots)$.
This automorphism
already appeared in the study of the
non-equivariant version of $K$-Peterson isomorphism in
\cite{BMS,IIN}.

\subsection[Image of g\^{}(k)\_la(y|b) under sigma for k-small la]{Image of $\boldsymbol{{g}^{(k)}_\la(y|b)}$ under $\boldsymbol{\sigma}$ for $\boldsymbol{k}$-small $\boldsymbol{\la}$}

Let $\mathrm{diag}(\lambda)$ denote the \textit{main diagonal} of $\lambda$, that is, the set of boxes at the $(i,i)$-th position for~${i=1,2,\dots}$.
For $x\in \mathrm{diag}(\lambda)$, let $r(x)$ be the $n$-residue of the box that is furthest to the right from $x$, and $b(x)$ the $n$-residue of the box that is furthest below from $x$.
\begin{Proposition}\label{prop:v1_of_kKShur}
If $\lambda$ is $k$-small, we have
\[
\sigma\bigl(g^{(k)}_\lambda(y|b)\bigr)
=
\Bigg(\prod_{x\in \mathrm{diag}(\lambda)}{\rm e}^{a_{r(x)+1}-a_{b(x)}}\Bigg)
\tilde{g}^{(k)}_\lambda(y|b).\]
\end{Proposition}

\begin{Example}
When $n=6$ and
\[
\ytableausetup{smalltableaux}\lambda=
 \begin{ytableau}
 0 & 1& 2\\ 5&0&1\\4
 \end{ytableau}
 ,
 \]
the main diagonal consists of two boxes $x_1$, $x_2$, where $x_i$ is at the $(i,i)$-th position.
Since
$r(x_1)=2$,
$b(x_1)=4$,
$r(x_2)=1$, and
$b(x_2)=0$,
we have \smash{$\sigma(g^{(5)}_{\lambda}(y|b))={\rm e}^{(a_3-a_4)+(a_2-a_6)}\Tilde{g}^{(5)}_\lambda(y|b)$}.
\end{Example}

\begin{Lemma}\label{lem:sigma T0}
$\sigma \circ T_0={\rm e}^{\theta} D_0 \circ \sigma$.
\end{Lemma}
\begin{proof}
Let \smash{$f(y|b)\in \hat{\Lambda}^\Rep$}. Then
by using~\eqref{eq:iota Omega}, we have
\begin{align*}
\sigma T_0 f(y|b)
&=\sigma
\left(
\frac{1}{1-{\rm e}^{-\theta}}(\frac{\Omega(b_1|y)}{\Omega(b_n|y)}s_\theta f(y|b) -f(y|b))
\right)\\
&=
\frac{1}{1-{\rm e}^{-\theta}}\left({\rm e}^{\theta}\frac{\Omega(b_1|y)}{\Omega(b_n|y)} \sigma(( s_\theta f)(y| b)) -(\sigma f)(y|b)\right)\\
&=
\frac{1}{1-{\rm e}^{-\theta}}\left({\rm e}^{\theta}\frac{\Omega(b_1|y)}{\Omega(b_n|y)}s_\theta(( \sigma f)(y| b))-{\rm e}^{\theta} \sigma f(y| b) +{\rm e}^{\theta}\sigma f(y| b)-\sigma f(y|b)\right)\\
&={\rm e}^{\theta}
T_0\sigma f(y|b)
+{\rm e}^{\theta}\sigma f(y|b)={\rm e}^{\theta}
D_0\sigma f(y|b),
\end{align*}
where we use that $\sigma$ commutes with the action of $s_\theta$.
\end{proof}

\begin{proof}[Proof of Proposition~\ref{prop:v1_of_kKShur}]
We use induction on the number of boxes of $\la$. The case $\la=\varnothing$ is
obvious. Suppose~${\la\ne \varnothing}$.
There is a~box removable from $\la$, with the $n$-residue say $i$. Let $\mu$ be the partition obtained from $\la$ by removing the box.
We consider the case when $i=0$.
One easily sees that~${r(x)+1}$, $b(x)\notin\{1,n\}$ for $x\in \mathrm{diag}(\mu)$. It follows that $s_\theta (e(\mu))=e(\mu)$. We also note that~${{\rm e}^\theta e(\mu)=e(\lambda)}$,
\begin{align*}
\sigma\bigl(g^{(k)}_\lambda(y|b)\bigr)
&=\sigma\bigl(T_0 g^{(k)}_\mu(y|b)\bigr)={\rm e}^{\theta}{D}_0\sigma\bigl( g^{(k)}_\mu(y|b)\bigr)\qquad \text{by Lemma~\ref{lem:sigma T0}}\\
&={\rm e}^{\theta}{D}_0\bigl(e(\mu)\tilde{g}^{(k)}_\mu(y|b)\bigr)={\rm e}^{\theta}e(\mu){D}_0\bigl(\tilde{g}^{(k)}_\mu(y|b)\bigr)
=e(\lambda)\tilde{g}^{(k)}_\lambda(y|b).
\end{align*}
The case when $i\ne 0$ is
left to the reader since it is
similar and easier.
\end{proof}

\begin{Remark}
The assumption that $\lambda$ is $k$-small in Proposition~\ref{prop:v1_of_kKShur} is mandatory.
Indeed, if we consider the case when $n=3$ and $w=s_2s_1s_0$, the associated $2$-bounded partition $\lambda=(2,1)$ is not $2$-small, and
\smash{$\sigma\bigl(g^{(2)}_w(y|b)\bigr)$} equals
\begin{gather*}
\mathcal{T}_2\bigl({\rm e}^{a_2-a_3}\otimes \tilde{g}^{(2)}_{10}(y|b)\bigr)\\
\qquad
=\bigl(1+{\rm e}^{a_3-a_2}\bigr)\otimes \tilde{g}^{(2)}_{10}(y|b)+{\rm e}^{a_3-a_2}\otimes T_2\tilde{g}^{(2)}_{10}(y|b) \\
\qquad=1\otimes \tilde{g}^{(2)}_{10}(y|b)+{\rm e}^{a_3-a_2}\otimes D_2\tilde{g}^{(2)}_{10}(y|b)
=1\otimes \tilde{g}^{(2)}_{10}(y|b)+{\rm e}^{a_3-a_2}\otimes \tilde{g}^{(2)}_{210}(y|b).
\end{gather*}
\end{Remark}

\subsection[Regular functions c\_i]{Regular functions $\boldsymbol{c_i}$}\label{sec:ci}

For arbitrary matrix $Z$ that commutes with $A$, there are unique scalars $c_i$ satisfying
\begin{equation}\label{eq:def_of_c_i}
Z=c_0E+c_1A+\cdots+c_{n-1}A^{n-1}.
\end{equation} The existence of such
scalars is assured by the fact that $A$ is conjugate to a~companion matrix.
We will consider each $c_i$ as an element of $\mathscr{O}\bigl(\tilde{\mathscr{Z}}\bigr)$.

\begin{Example} If
$n=3$ we have expressions
\[
c_0=z_{11}+{\rm e}^{-a_1}z_{12}+{\rm e}^{-a_1-a_2}z_{13},\qquad
c_1=-z_{12}-({\rm e}^{-a_1}+{\rm e}^{-a_2})z_{13},\qquad
c_2=z_{13}
\]
by comparing the $1$st row of~\eqref{eq:def_of_c_i}.
\end{Example}

By comparing the diagonal
entries, we have
\smash{$
z_{ii}=c_0+c_1 {\rm e}^{-a_i}+
\cdots+c_{n-1}{\rm e}^{-(n-1)a_i}$}.
We also define \smash{$c_0^{(j)},c_1^{(j)},\dots,c_{n-1}^{(j)}\in \mathscr{O}\bigl(\tilde{\mathscr{Z}}\bigr)$}
for $j\in \Z$ by
\begin{equation}\label{eq:def_of_c_i_(j)}
A^jZ=c_0^{(j)}E+c_1^{(j)}A+\cdots+c_{n-1}^{(j)}A^{n-1}.
\end{equation}
Note that we have $c_j^{(0)}=c_j$.
Comparing the diagonal entries on both sides of~\eqref{eq:def_of_c_i_(j)}, we have
\begin{equation}
{\rm e}^{-ja_i}z_{ii}=c_0^{(j)}+c_1^{(j)} {\rm e}^{-a_i}+
\cdots+c_{n-1}^{(j)}{\rm e}^{-(n-1)a_i}.\label{eq:zii_c_2}
\end{equation}

\begin{Theorem}\label{thm:det_c_to_tausigma}
We have
\begin{equation}\label{eq:det_c_to_tausigma}
\begin{gathered}
\sigma_i=
\begin{vmatrix}
c^{(0)}_{0} & c^{(0)}_{1} & \dots & c^{(0)}_{i-1}\\
c^{(1)}_{0} & c^{(1)}_{1} & \dots & c^{(1)}_{i-1}\\
\vdots & \vdots & \dots & \vdots \\
c^{(i-1)}_{0} & c^{(i-1)}_{1} & \dots & c^{(i-1)}_{i-1}\\
\end{vmatrix}, \qquad
\tau_i=
\begin{vmatrix}
c^{(1)}_{0} & c^{(1)}_{1} & \dots &c^{(1)}_{i-1}\\
c^{(2)}_{0} & c^{(2)}_{1} & \dots &c^{(2)}_{i-1}\\
\vdots & \vdots & \dots & \vdots\\
c^{(i)}_{0} & c^{(i)}_{1} & \dots &c^{(i)}_{i-1}\\
\end{vmatrix}.
\end{gathered}
\end{equation}
\end{Theorem}

\begin{proof}
Since $P^{-1}$ is lower unitriangular, we use
Proposition~\ref{prop:Noumi} to have
\[
\begin{aligned}
\sigma_i=\det(ZP)^{[1,i]}_{[1,i]}
=
\det\bigl(P^{-1}ZP\bigr)^{[1,i]}_{[1,i]}=
\det\bigl(c_0E+c_1C_A+\dots+c_{n-1}C_A^{n-1}\bigr)^{[1,i]}_{[1,i]}.
\end{aligned}
\]
By comparing the 1-st row of
\[
C_A^{j-1}\bigl(c_0E+c_1C_A+\cdots+c_{n-1}C_A^{n-1}\bigr)
=c_0^{(j-1)}E+c_1^{(j-1)}C_A+\cdots+c_{n-1}^{(j-1)}C_A^{n-1},
\]
we see that the
$(j,l)$-th entry of $c_0E+c_1C_A+\cdots+c_{n-1}C_A^{n-1}$
is
\smash{$(-1)^{j+l}c_l^{(j-1)}$}.
Hence the first equation of~\eqref{eq:det_c_to_tausigma} holds.
The equation for $\tau_i$ follows from this by applying $\sigma^{-1}$ because of~\eqref{eq:sigma_tau_i}.\looseness=-1
\end{proof}

\subsection[Involution iota on O(Z)]{Involution $\boldsymbol{\iota}$ on $\boldsymbol{\mathscr{O}\bigl(\tilde{\mathscr{Z}}\bigr)}$}
\label{sec:iota on OZ}

Let $J$ be the permutation matrix
of the longest element $w_\circ$ of $S_n$.
Explicitly, $J$ is $\sum_{i=1}^n E_{i,n-i+1}$.
Let $\iota$ be a~ring automorphism of $\mathscr{O}\bigl(\tilde{\mathscr{Z}}\bigr)$ defined by
$\iota({\rm e}^\gamma)=
{\rm e}^{-w_\circ \gamma}$ (in particular
$\iota({\rm e}^{a_i})={\rm e}^{-a_{n-i+1}}$), and
\[
\sum_{i=0}^{n-1}
\iota(c_i)A^{-i}
=\left(\sum_{i=0}^{n-1}
c_iA^{i}\right)^{-1}.
\]
It is not difficult to show $\iota(C_A)=JC_A^{-1}J$.
Since $C_A=P^{-1}AP$, we obtain
\[\sum_{i=0}^{n-1}
\iota(c_i)C_A^{-i}
=\left(\sum_{i=0}^{n-1}
c_iC_A^{i}\right)^{-1}
\] and
\begin{align}
\iota\bigl(P^{-1}ZP\bigr)
&=\iota\bigl(c_0+c_1C_A+\dots+c_{n-1}C_A^{n-1}\bigr) \nonumber\\
&=J\bigl(\iota(c_0)+\iota(c_1)C_A^{-1}+\dots+\iota(c_{n-1})C_A^{-(n-1)} \bigr)J\nonumber\\
&=J\bigl(c_0+c_1C_A+\dots+c_{n-1}C_A^{n-1} \bigr)^{-1}J=JP^{-1}Z^{-1}PJ.\label{eq:S_PZP}
\end{align}
Comparing the diagonal entries on both sides of~\eqref{eq:S_PZP}, we have
\begin{equation}
\label{eq:S_zii}
\iota(z_{ii})=z_{n-i+1,n-i+1}^{-1}.
\end{equation}

Identifying $z_{ii}$ with $\Omega(b_i|y)$, we have
\begin{equation}\label{eq:iota Omega}
\iota(\Omega(b_i|y))=\Omega(b_{n-i+1}|y)^{-1}.
\end{equation}

\begin{Remark}
The automorphism $\iota$ is induced by a~diagram automorphism of the affine Dynkin diagram. See \cite[Sections~2.7 and~3.3]{ISY}.
\end{Remark}

\begin{Proposition}\label{prop:S_and_sigma}
$\iota$ commutes with $\sigma$.
\end{Proposition}
\begin{proof}
Note that $\iota$ and $\sigma$ are naturally extended
as ring automorphisms
of $\mathscr{O}\bigl(\tilde{\mathscr{Z}}\bigr)^\Delta$. In view of
Proposition~\ref{prop:OZ_Delta}, we only need to check on the generators $z_{ii}$
and $f\in R(T)$.
From~\eqref{eq:S_zii} and~\eqref{eq:sigma_zii}, we have
\begin{align*}
 \sigma (\iota(z_{ii}))
 =\sigma\bigl(z_{n-i+1,n-i+1}^{-1}\bigr)=
 {\rm e}^{-a_{n-i+1}}
 z_{n-i+1,n-i+1}^{-1}=\iota({\rm e}^{a_i}z_{ii})=\iota(\sigma(z_{ii})).
\end{align*}
Since $\sigma$ is $R(T)$-linear, it is clear that
$\iota(\sigma(f))=\sigma(\iota(f))=\iota(f)
$ for $f\in R(T)$.
\end{proof}

\begin{Proposition}
\label{prop:S_on_tau}
We have
\begin{align}
\label{eq:sigma tau_i}
\iota(\tau_i)&=\frac{\tau_{n-i}}{\tau_n},\qquad
\iota(\sigma_i)=\frac{\sigma_{n-i}}{\sigma_n}.
\end{align}
\end{Proposition}
\begin{proof}
By~\eqref{eq:S_PZP}, we have
\begin{align*}
\iota(\sigma_i)
&=\iota\left(\det\bigl(P^{-1}ZP\bigr)_{[1,i]}^{[1,i]}\right)=\det\bigl(JP^{-1}Z^{-1}PJ\bigr)_{[1,i]}^{[1,i]}=\det\bigl(P^{-1}Z^{-1}P\bigr)_{[n-i+1,n]}^{[n-i+1,n]}\\
&=
\det\bigl(P^{-1}Z^{-1}P\bigr)\cdot\det\bigl(P^{-1}ZP\bigr)_{[1,n-i]}^{[1,n-i]}=
\det(Z^{-1})\cdot\det\bigl(P^{-1}ZP\bigr)_{[1,n-i]}^{[1,n-i]}\\
&=\frac{\sigma_{n-i}}{\Omega(b_1|y)\cdots \Omega(b_n|y)}
=\frac{\sigma_{n-i}}{\sigma_n},
\end{align*}
where for the fourth equality
we use
the fact
\smash{$\det A^{[m+1,n]}_{[m+1,n]}
=\det A \cdot\det\bigl(A^{-1}\bigr)^{[1,m]}_{[1,m]}$},
which holds for any invertible $A$.
The first equality of~\eqref{eq:sigma tau_i} is obtained from Proposition~\ref{prop:S_and_sigma} and~\eqref{eq:sigma_tau_i}.
\end{proof}

\subsection[Basic properties of Phi\_n]{Basic properties of $\boldsymbol{\Phi_n}$}
We collect some basic properties of $\Phi_n$ which will be used below.

\begin{Proposition}\label{prop:S_commutes}
We have
\begin{align}
 \iota\circ \Phi_n&=\Phi_n \circ \iota.\label{eq:S_Phi}
\end{align}
\end{Proposition}
\begin{proof}
 For~\eqref{eq:S_Phi}, it suffices to show on generators $z_i$, and $Q_i$.
From Proposition~\ref{prop:S_on_tau}, we have
\[
(\iota\circ \Phi_n)(z_i)
=\iota\left(\frac{\tau_i\sigma_{i-1}}{\tau_{i-1}\sigma_i}\right)
=\frac{\tau_{n-i}\sigma_{n-i+1}}{\sigma_{n-i}\tau_{n-i+1}}=\Phi_n(z_{n-i+1}^{-1})
=(\Phi_n\circ \iota)(z_i)\]
and
\[
(\iota\circ\Phi_n)(Q_i)=\iota\left(\frac{\tau_{i-1}\tau_{i+1}}{\tau^2_{i}}\right)=
\frac{\tau_{n-i+1}\tau_{n-i-1}}{\tau^2_{n-i}}=\Phi_n(Q_{n-i})
=(\Phi_n\circ \iota)(Q_i),
\]
which concludes the proposition.
\end{proof}

Let us denote
the characteristic polynomial of \smash{$L_{[1,i]}^{[1,i]}$}
by
\smash{$
\chi_i(\zeta):=\det(\zeta E-L)_{[1,i]}^{[1,i]}$}.
Then we~have from Appendix~\ref{sec:F}
\[
\chi_i(\zeta)=
\zeta^i-F^{(i)}_1\zeta^{i-1}+\dots+(-1)^{i}F^{(i)}_i,
\]
where
\begin{equation}
\label{eq:Fij}
F_m^{(i)}:=\sum_{\substack{J\subset [i]\\ |J|=m}}
\prod_{\substack{j\in J,\\j+1\notin J}}
(1-Q_j)\prod_{j\in J}z_j,
\end{equation}
for
$0\le m\le i\le n$, with $F_0^{(i)}=1$.

\begin{Lemma}\label{lem:SF_to_c}
For $1\le m\le i\le n$, we have
\[
(\iota\circ\Phi_n )\bigl(F^{(i)}_m\bigr)=\frac{(-1)^{m}}{\tau_{n-i}}
\begin{vmatrix}
c^{(1)}_{0} & c^{(1)}_{1} & \dots & c^{(1)}_{n-i-2}&c^{(1)}_{n-i+m-1}\\
c^{(2)}_{0} & c^{(2)}_{1} & \dots & c^{(2)}_{n-i-2}&c^{(2)}_{n-i+m-1}\\
\vdots & \vdots & \dots & \vdots &\vdots\\
c^{(n-i)}_{0} & c^{(n-i)}_{1} & \dots & c^{(n-i)}_{n-i-2}&c^{(n-i)}_{n-i+m-1}\\
\end{vmatrix}.
\]
\end{Lemma}
 \begin{proof}
We first compute $\Phi_n(\chi_i(\zeta))$ as follows
\begin{gather*}
\det\bigl(U^{-1}\bigr)_{[1,i]}^{[1,i]}\cdot
\Phi_n(\chi_i(\zeta))
\cdot\det R_{[1,i]}^{[1,i]}\\
\qquad=
\det\bigl(U^{-1}\bigr)_{[1,i]}^{[1,i]}\cdot
\det(\zeta E-\Phi_n(L))_{[1,i]}^{[1,i]}
\cdot\det R_{[1,i]}^{[1,i]}\\
\qquad=
\det\bigl( U^{-1}(\zeta E-\Phi_n(L))R\bigr)_{[1,i]}^{[1,i]}\qquad\text{by Proposition~\ref{prop:Noumi}}
\\
\qquad=
\det\bigl((\zeta E-C_A) U^{-1}R\bigr)_{[1,i]}^{[1,i]}
\qquad\text{since $\Phi_n(L)=UC_AU^{-1}$ from
\eqref{eq:AZ_to_L}}\\
\qquad={\det\bigl((\zeta E-C_A)P^{-1}ZAP\bigr)_{[1,i]}^{[1,i]}} \qquad\text{by~\eqref{eq:Gauss_decomp}}.
\end{gather*}
Since \smash{$\det\bigl(U^{-1}\bigr)_{[1,i]}^{[1,i]}=1$} and
\smash{$\det R_{[1,i]}^{[1,i]}=\tau_i$} (see~\eqref{eq:r_to_minor}), we deduce that
\begin{equation}\label{eq:PhiL_Z}
\Phi_n(\chi_i(\zeta))
={\tau_i}^{-1}{\det\bigl((\zeta E-C_A)P^{-1}ZAP\bigr)_{[1,i]}^{[1,i]}}.
\end{equation}
Let us apply $\iota$ on both sides of~\eqref{eq:PhiL_Z}.
By using Proposition~\ref{prop:S_commutes}, we can verify
\[
\iota\bigl(P^{-1}ZAP\bigr)=
JP^{-1}Z^{-1}A^{-1}PJ.
\]
Using this together with $\iota(C_A)=JC_A^{-1}J$, we have
\begin{align*}
(\iota\circ \Phi_n)(\chi_i(\zeta))
&=(\tau_n/\tau_{n-i})\cdot{\det\bigl(J(\zeta E-C_A^{-1})P^{-1}Z^{-1}A^{-1}PJ\bigr)_{[1,i]}^{[1,i]}}\\
&=(\tau_n/\tau_{n-i})\cdot\det\bigl(\bigl(\zeta E-C_A^{-1}\bigr)P^{-1}Z^{-1}A^{-1}P\bigr)_{[n-i+1,n]}^{[n-i+1,n]}\\
&=\frac{1}{\tau_{n-i}}
\det\left(
\begin{matrix}
\bigl(P^{-1}AZP\bigr)_{[1,n-i]}^{[1,n]}\\
\hline
\bigl(\zeta E-C_A^{-1}\bigr)_{[n-i+1,n]}^{[1,n]}
\end{matrix}
\right),
\end{align*}
where in the
last equality we use Proposition~\ref{prop:det AB partial}, and $\det\bigl(P^{-1}AZP\bigr)=\tau_n$.
Comparing the coefficients of $\zeta^{i-m}$ on both sides of this equality, we obtain
\[
(\iota\circ\Phi_n)\bigl(F^{(i)}_m\bigr)=
\frac{\det\bigl(P^{-1}AZP\bigr)_{[1,n-i]}^{[1,n-i-1]\cup\{n-i+m\}}}{\tau_{n-i}}.
\]
Because we know $P^{-1}AZP=\bigl((-1)^{j+l}c^{(j)}_{l}\bigr)_{1\le j,l\le n}$ (see the proof of Theorem~\ref{thm:det_c_to_tausigma}), we obtain the desired result, after eliminating unnecessary signs.
\end{proof}

\section{Quantum double Grothendieck polynomials}\label{sec:Groth}
We use notation of binary operations $x\oplus y=x+y-xy$ and
$x\ominus y=(x-y)/(1-y)$.
We also denote $-x/(1-x)=0\ominus x$ by $\overline{x}$.

\subsection{Quantum double Grothendieck polynomials}
\label{sec:G^Q(z|eta)}
We recall
Lenart--Maeno's
quantum double Grothendieck polynomials.

Let us consider
the two sets of variables
$z=(z_1,\dots,z_n)$ and
$\eta=(\eta_1,\dots,\eta_n)$, and the polynomial ring
$\Z[Q][z_1,\dots,z_n,\eta_1,\dots,\eta_n]$,
where $Q_1,\dots,Q_{n-1}$ are the Novikov variables.
For $i\in I$, define
\begin{equation}
\label{eq:DiQ,TiQ}
T_i^Q
=\frac{s_i^{(\eta)}-1}{\eta_{i+1}\ominus \eta_{i}},\qquad
D_i^Q
=1+\frac{s_i^{(\eta)}-1}{\eta_{i+1}\ominus \eta_{i}},
\end{equation}
as linear endomorphisms of
$\Z[Q][z_1,\dots,z_n,\eta_1,\dots,\eta_n]$.
Define
\begin{gather}
\label{eq:G_longest}
\Gro{\wz}(z|\eta)=
\prod_{i=1}^{n-1}
\fac_i,\\
\fac_i=\sum_{j=0}^i
(-1)^j
(1-\eta_{n-i})^j
F_j^{(i)}(z,Q).
\label{eq:def phi}
\end{gather}
Note that $\fac_i$ can be written as
\begin{equation}
\label{eq:phi det}
\fac_i=\det(E-(1-\eta_{n-i})L)^{[1,i]}_{[1,i]}
\end{equation}
if $L$ is the matrix given as in Section~\ref{sec:Rel Toda}.

There exists a~unique collection of polynomials
\[
\bigl\{\Gro{w}(z|\eta)
\in \Z[Q][z_1,\dots,z_{n},\eta_1,\dots,\eta_{n}]
 \mid w\in S_n\bigr\},\]
where $\mathfrak{G}_{w_\circ}^Q$ corresponding to the longest element $w_\circ$ is given by~\eqref{eq:G_longest}, and each $\mathfrak{G}_{w}^Q$ for general~$w$ is characterized by the recursive relation
\begin{equation}
 \label{eq:Di_Gr}
D_i^Q\Gro{w}(z|\eta)=\begin{cases}
\Gro{s_iw}(z|\eta) & \text{if $s_i w<w$},\\
\Gro{w}(z|\eta) & \text{if $s_i w>w$}.
\end{cases}
\end{equation}
\begin{Remark}\label{rem:e_vs_eta}
We follow the notation in~\cite{MNS2} as closely as possible, however, there are some
unavoidable differences as outlined below.
The variables $y_i$ used for equivariant parameters in~\cite{MNS2} are denoted here by $\eta_i$.
This change is necessary because we reserve $y_i$ for the variables of symmetric functions.
We treat the variable $\eta_i$ as
an element of $\Rep$ via the correspondence
\begin{equation}
\eta_i=1-{\rm e}^{a_{n-i+1}}=\overline{b}_{n-i+1},\qquad1\le i\le n.
\label{eq:eta-b}
\end{equation}
Note that our $a_i$ corresponds to $-\epsilon_i$ in~\cite{MNS2}.
Additionally, we use $z_i=1-x_i$ with $x_i$
the variable from~\cite{MNS2}.

\end{Remark}

\begin{Proposition}For $i\in I$, we have
 \begin{equation}
D_{n-i}\circ \Phi_n=\Phi_n \circ D_i^Q,\qquad
D_{n-i}\circ \tilde{\Phi}_n=\tilde{\Phi}_n \circ D_i^Q.
\label{eq:D_and_Phi}
\end{equation}
\end{Proposition}
\begin{proof}
The operators on both sides
coincide on $z_i$, $Q_i$,
because $\tau_i$ and $\sigma_i$ are $S_n$-invariant. Since~$\Phi_n$ \big(and $\tilde{\Phi}_n$\big) is $\Rep$-linear,
it suffices to see the equality
on $\Rep$. In fact, we have
\begin{equation*}
1+\frac{s_i^{(\eta)}-1}{\eta_{i+1}\ominus \eta_{i}}
=1+
\frac{s_{n-i}^{(b)}-1}{\overline{b}_{n-i}\ominus \overline{b}_{n-i+1}}
=1+
\frac{s_{n-i}^{(b)}-1}{b_{n-i+1}\ominus b_{n-i}},
\end{equation*}
under the
identification~\eqref{eq:eta-b},
where \smash{$s_i^{(b)}$} exchanges $b_i$
and $b_{i+1}$,
and we use $\overline{x}\ominus \overline{y}=y\ominus x$.
\end{proof}

\subsection[Involution iota on the equivariant quantum K-ring]{Involution $\boldsymbol{\iota}$ on the equivariant quantum $\boldsymbol{K}$-ring}

Let $\iota$ be a~ring homomorphism of
$\Rep[Q]\big[z_1^{\pm 1},\dots,z_n^{\pm 1}\big]$
such that
\[
\iota({\rm e}^{a_i})={\rm e}^{-a_{n-i+1}},\qquad
\iota(z_{i})=z_{n-i+1}^{-1},\qquad
\iota(Q_i)=Q_{n-i}.
\]
Obviously,
$\iota$ satisfies $\iota^2=1$.
It is easy to prove that
\begin{equation}\label{eq:S_and_D_i}
 \iota\circ D_i^Q=D_{n-i}^Q\circ \iota.
 \end{equation}
 holds for $i\in I$.
\begin{Proposition}
$\iota$ preserves the ideal $\mathscr{I}_n^Q$ of $\Rep[Q]\big[z_1^{\pm 1},\dots,z_n^{\pm 1}\big]$.
\end{Proposition}
\begin{proof}
We claim that
\smash{$
z_1\cdots z_n \iota\bigl(F_i^{(n)}\bigr)=F_{n-i}^{(n)}$}, \smash{$
e_n^{(n)}\iota\bigl(e_i^{(n)}\bigr)=e_{n-i}^{(n)}
$}
for $1\le i\le n-1$.
The second identity is obvious.
For $J\subset [1,n]$, $ |J|=i$, set
$K:=\{i'\mid i\notin J\}$, with $i':=n-i+1$. Then
$z_1\cdots z_n \prod_{j\in J}\iota(z_i)=\prod_{l\in K}z_l$.
We see that $j\in J, j+1\notin J$
is equivalent to
$(j+1)'\in K$, $j'\notin K$.
Thus
\[
\prod_{j\in J, j+1\notin J}
(1-Q_j)=
\prod_{l\in K, l+1\notin K}
(1-Q_l).
\]
Hence the claim holds and we have $z_1\cdots z_n \iota\bigl(F_i^{(n)}\bigr)=F_{n-i}^{(n)}$.
Since $F_n^{(n)}=z_1\cdots z_n$, and $e_n^{(n)}={\rm e}^{-(a_1+\cdots+a_n)}=1$, we are done.
\end{proof}

The next result will be proved in Section~\ref{sec:max factor}.
\begin{Proposition}\label{prop:long_iota_invariant}
We have $\iota\bigl(\mathfrak{G}_{w_\circ}^{Q}\bigr)=\mathfrak{G}_{w_\circ}^{Q}$.
\end{Proposition}

For $w\in S_n$, let
$w^*=w_\circ w w_\circ;$
$w^*$ is obtained from any
reduced expression of $w$ by replacing~$s_i$ with $s_{n-i}$.
Note that the map $w\mapsto w^*$ preserves the Bruhat order.
\begin{Proposition}\label{prop:S_star}
Let $w\in S_n$. Then
\begin{equation}
\iota\bigl(\Gro{w}(z|\eta)\bigr)=
\Gro{w^*}(z|\eta).
\label{eq:S_star}
\end{equation}
\end{Proposition}
\begin{proof}
We use induction on $\ell(w_\circ)-\ell(w)$.
Since $\Gro{w_\circ}$
is $\iota$-invariant (see Proposition~\ref{prop:long_iota_invariant})
and $w_\circ^*=w_\circ$, \eqref{eq:S_star} holds for $w=w_\circ$.
Suppose $w\in S_n$ satisfies
$s_iw>w$ for some $1\le i\le n-1$.
Then we have $(s_iw)^*>w^*$.
It follows that $s_{n-i}(s_iw)^*<(s_iw)^*$.
Then by using
\eqref{eq:Di_Gr} and
\eqref{eq:S_and_D_i}, we have
\begin{align*}
\iota\bigl(\Gro{w}(z|\eta)\bigr)&=
\iota\bigl(D_i^Q \Gro{s_iw}(z|\eta)\bigr)=D_{n-i}^Q \iota\bigl(\Gro{s_iw}(z|\eta)\bigr)\\
&=D_{n-i}^Q \Gro{(s_iw)^*}(z|\eta)=\Gro{s_{n-i}(s_iw)^*}(z|\eta)=\Gro{w^*}(z|\eta).\tag*{\qed}
\end{align*}\renewcommand{\qed}{}
\end{proof}

\subsection[Explicit formula for Gro\{s\_theta\}(x|eta)]{Explicit formula for $\boldsymbol{\Gro{s_\theta}(x|\eta)}$}
The quantum double Grothendieck polynomial for $s_\theta$ plays a~fundamental role
in the following.

Let $\omega^{(\eta)}$ be the cyclic permutation with respect to the variables $\eta$, that is
$\omega^{(\eta)}(\eta_i)=\eta_{i+1}$
with $\eta_{n+1}=\eta_1$.
We use the following notation:
$
[x|\eta]^i=(x\oplus \eta_1)\cdots (x\oplus \eta_i)
$
for $i\ge 1$.
\begin{Proposition}
\label{prop:Gtheta}
We have
\begin{equation}
\Gro{s_\theta}(z|\eta)=\Gro{s_1s_2\dots s_{n-1}}(z|\eta)\cdot \Gro{s_{n-2}\dots s_2s_1}\bigl(z|\omega^{(\eta)}(\eta)\bigr)\label{eq:Gtheta}
.
\end{equation}
\end{Proposition}
\begin{proof}
We first prove
the classical version
\begin{equation}
\mathfrak{G}_{s_\theta}(z|\eta)=\mathfrak{G}_{s_1s_2\dots s_{n-1}}(z|\eta)\cdot \mathfrak{G}_{s_{n-2}\dots s_2s_1}\bigl(x|\omega^{(\eta)}(\eta)\bigr)
\label{eq:Gtheta-cl}
\end{equation}
of
\eqref{eq:Gtheta}.
Since
$s_\theta\in S_n$ is a~dominant permutation of
code (see \cite[Chapter~I]{Mac} for the definitions) $c(s_\theta)=(n-1,1,\dots,1,0)$,
\begin{align*}
\mathfrak{G}_{s_\theta}(z|\eta)&=
[x_1|\eta]^{n-1}
[x_2|\eta]
\cdots
[x_{n-1}| \eta]=[x_1|\eta]
[x_2|\eta]
\cdots
[x_{n-1}| \eta]\times
\big[x_1|\omega^{(\eta)}(\eta)\big]^{n-2}.
\end{align*}
For the first equality, we have used \cite[Lemma B.6]{MNS1}.
Because
$s_1s_2\dots s_{n-1}$
and $s_{n-2}\dots s_{2}s_1$
are also dominant, with the codes
$(1,\dots,1,0)$ and $ (n-2,0,\dots,0)$
respectively, we have by \cite[Lemma B.6]{MNS1}
\begin{gather*}
\mathfrak{G}_{s_1s_2\dots s_{n-1}}(z|\eta)
=[x_1|\eta]
[x_2|\eta]
\cdots
[x_{n-1}| \eta],\qquad
\mathfrak{G}_{s_{n-2}\dots s_{2}s_1}(z|\eta)
=[x_1|\eta]^{n-2}.
\end{gather*}
Therefore, \eqref{eq:Gtheta-cl} holds.

We will apply the quantization map
$\hat{Q}$ of Lenart--Maeno (see \cite[Section~3]{LM}) to~\eqref{eq:Gtheta-cl} with respect to the variables
$x_1,\dots,x_n$, where $x_i=1-z_i$.
Let \smash{$f_j^{(i)}=e_j(1-x_1,\dots,1-x_i)$}. Then
$\hat{Q}$ is a~linear map such that
\[\hat{Q}\bigl(f_{p_1}^{(1)}f_{p_2}^{(2)}\cdots f_{p_{n-2}}^{(n-2)}f_{p_{n-1}}^{(n-1)}\bigr)=F_{p_1}^{(1)}F_{p_2}^{(2)}\cdots F_{p_{n-2}}^{(n-2)}F_{p_{n-1}}^{(n-1)}
\]
for $0\le p_i\le i$ (see \cite[Proposition 3.16]{LM}).
As a~polynomial in $x_1$, $[x_1|\eta]^{n-2}$
has degree~${n-2}$.
Hence it is easy to see that
$\mathfrak{G}_{s_{n-2}\dots s_{2}s_1}(z|\eta)=[x_1|\eta]^{n-2}$
is a~linear combination of
$\smash{f_{p_1}^{(1)}f_{p_2}^{(2)}}\cdots\allowbreak\smash{ f_{p_{n-2}}^{(n-2)}}$ with $ 0\le p_i\le i$.
Also,
\begin{equation}
\mathfrak{G}_{s_{1}s_2\dots s_{n-1}}(z|\eta)=[x_1|\eta]
[x_2|\eta]
\cdots
[x_{n-1}| \eta]=
1+\sum_{j=1}^{n-1}
(\eta_1-1)^jf_j^{(n-1)}.
\label{eq:phi_(n-1)-2}
\end{equation}
Therefore,
$\hat{Q}$ preserves the product of~\eqref{eq:Gtheta-cl},
so we have~\eqref{eq:Gtheta}.
\end{proof}

By the proof of the previous
proposition, we have the following.
\begin{Corollary}\label{cor:phi_(n-1)}
We have
\[
\Gro{s_1s_2\dots s_{n-1}}(z|\eta)=
\sum_{j=0}^{n-1}(-1)^j(1-\eta_{1})^{j}F_j^{(n-1)}.
\]
\end{Corollary}
\begin{proof} This identity is obtained by
applying $\hat{Q}$ to
\eqref{eq:phi_(n-1)-2}.
\end{proof}

\subsection[Affine K-nil-Hecke action on the equivariant quantum K-ring]{Affine $\boldsymbol{K}$-nil-Hecke action on the equivariant quantum $\boldsymbol{K}$-ring}
Here we explain an action of
the $0$-th Demazure operator on
$QK_T(\SL_n(\C)/B)_Q$.
We define an operator \smash{$D_0^Q$} on
$QK_T(\SL_n(\C)/B)_Q$ by
\begin{equation}
\label{eq:D0Q}
D_0^Q=T^{Q}_{\theta}+
Q^{-\theta^\vee}
\mathscr{O}^{s_\theta}\cdot
s_\theta,
\end{equation}
where an element $T^Q_\theta$ is defined by
\[
T^Q_\theta:=\bigl(\omega^{(\eta)}\bigr)^{-i}\circ T^Q_i\circ \bigl(\omega^{(\eta)}\bigr)^i=\frac{s_\theta^{(\eta)}-1}{\eta_1\ominus\eta_n}=\frac{s^{(\eta)}_\theta-1}{1-{\rm e}^{-\theta}}.
\]

If we present $QK_T(\SL_n(\C)/B)$ as a~quotient ring by~\cite{MNS1}, we can replace
$Q^{-\theta^\vee}\mathscr{O}^{s_\theta}$, acting as a~multiplication operator,
with \smash{$Q^{-\theta^\vee}\Gro{s_\theta}$}.

Together with $D_i^Q$, $i\in \Iaf$, and the left multiplication by $R(T)$,
 $QK_T(\SL_n(\C)/B)_Q$ has a~structure of $K$-nil-DAHA-module (see~\cite{KNS,KoNOS,Orr}).

\begin{Proposition}\label{prop:Quantum_D(1)}
Let $x=wt_\xi\in \hat{W}_G^0$. Then, with the presentation~{\rm\cite{MNS1}}, we have
\[
D_x^Q(1)
=Q^\xi \Gro{w} \mod \mathscr{I}_n^Q.
\]
\end{Proposition}
See Appendix~\ref{sec:nil-Hecke_on_QK}.

\subsection{Proof of Theorem~\ref{thm:main}}
The purpose of this section is to prove Theorem~\ref{thm:main}
\begin{Remark}
The image of $\tilde{\Phi}_n$ is \smash{$\hat{\Lambda}_{(n)}^\Rep\big[\sigma_i^{-1},(\sigma(\sigma_i))^{-1}|1\leq i\leq n\big]$}.
\end{Remark}

The following result is an important special case of Theorem~\ref{thm:main}.
This proposition in turn implies the commutativity of $\Phi_n$ with the Demazure operators (see Corollary~\ref{cor:Phi_commutes_w_Dem}).

\begin{Proposition}
\label{prop:key} We have
\[
\tilde{\Phi}_n\bigl(Q^{-\theta^\vee}\Gro{s_\theta}\bigr)=\tilde{g}_{s_0}^{(k)}(y|b).
\]
\end{Proposition}
\begin{Remark}The proof of this result is purely combinatorial.
We do not have to resort to
the deep geometric fact of Corollary~\ref{cor:BF}.
\end{Remark}

We first compute
the image of
the factor
$\Gro{s_1s_2\dots s_{n-1}}$ of $\Gro{s_\theta}$
in the factorized form~\eqref{eq:Gtheta} under~$\Phi_n$.
\begin{Lemma}\label{lem:Phi_1st_factor}
We have
\begin{gather}
 (\Phi_n\circ \iota)\bigl(\Gro{s_1s_2\dots s_{n-1}}\bigr)=
\frac{{\rm e}^{-a_1}z_{11}}{\tau_1},\label{eq:Phi_1st_factor}
 \\
 \Phi_n\bigl(\Gro{s_1s_2\dots s_{n-1}}\bigr)=
 \frac{\tau_n}{\tau_{n-1}}
 \frac{{\rm e}^{a_n}}{z_{nn}}.\label{eq:Phi_1st_factor-2}
\end{gather}
\end{Lemma}
\begin{proof}
 We will show~\eqref{eq:Phi_1st_factor}. Then
 \eqref{eq:Phi_1st_factor-2} is obtained from this by Proposition~\ref{prop:S_on_tau}.
Since $\iota(1-\eta_{1})={\rm e}^{-a_1}$, we have by Corollary~\ref{cor:phi_(n-1)}
\[
\iota\bigl(\Gro{s_1s_2\dots s_{n-1}}\bigr)=
\sum_{j=0}^{n-1}(-1)^j{\rm e}^{-ja_1}\iota\bigl(F_j^{(n-1)}\bigr).
\]
 By Lemma~\ref{lem:SF_to_c}, we have
\smash{$(\Phi_n\circ \iota)
\bigl(F^{(n-1)}_j\bigr)=(\iota\circ \Phi_n)
\bigl(F^{(n-1)}_j\bigr)=(-1)^j{c_j^{(1)}}/{\tau_1}$}.
Note that we have~\smash{$\tau_1=c_0^{(1)}$} from~\eqref{eq:det_c_to_tausigma}.
Therefore, we have
\begin{align*}
(\Phi_n\circ \iota)
\bigl(\Gro{s_1s_2\dots s_{n-1}}\bigr)
&=\frac{1}{\tau_1}
\sum_{j=0}^{n-1}{\rm e}^{-ja_1}
c_j^{(1)}=\frac{{\rm e}^{-a_1}z_{11}}{\tau_1}
\qquad \text{by~\eqref{eq:zii_c_2}}.\tag*{\qed}
\end{align*}\renewcommand{\qed}{}
\end{proof}
 In order to compute the image of
 $\Gro{s_{n-2}\dots s_2s_{1}}\bigl(z|\omega^{(\eta)}\eta\bigr)$ under $\Phi_n$, we prepare
 the following lemma.

\begin{Proposition}\label{prop:Phi_n_2nd_factor}
We have
\[
\Phi_n\bigl(\Gro{s_{n-2}\dots s_2s_{1}}(z|\eta)\bigr)=\frac{{\rm e}^{-a_1-a_2}z_{12}}{\tau_1}.
\]
\end{Proposition}
\begin{proof}
By~\eqref{eq:Di_Gr} and Proposition~\ref{prop:S_star}, we have
\begin{equation}
\Gro{s_{n-2}\dots s_2s_{1}}(z|\eta)
=
\iota\bigl(D_1^{Q}\Gro{s_{1}s_2\dots s_{n-1}}(z|\eta)\bigr).\label{eq:2nd_factor}
\end{equation}
Applying
${\Phi}_n$ to both sides of~\eqref{eq:2nd_factor}, we have
\begin{align*}
\Phi_n\bigl(\Gro{s_{n-1}\dots s_2s_{1}}(z|\eta)\bigr)
&=(\Phi_n\circ \iota)\bigl(D_1^{Q}\Gro{s_{1}s_2\dots s_{n-1}}(z|\eta)\bigr) \\
&=
D_1\bigl(
(\Phi_n\circ \iota)\bigl(\Gro{s_{1}s_2\dots s_{n-1}}(z|\eta)\bigr)
\bigr)\qquad\text{by~\eqref{eq:S_and_D_i} and~\eqref{eq:D_and_Phi}}\\
&=D_1\left(
\frac{{\rm e}^{-a_1}z_{11}}{\tau_1}
\right)\qquad \text{by
\eqref{eq:Phi_1st_factor}}\\
&=
\frac{1}{\tau_1}
D_1({{\rm e}^{-a_1}z_{11}})
\qquad \text{since $\tau_1$ is $S_n$-invariant.}
\end{align*}
Finally, it is easy to show
$D_1({\rm e}^{-a_1}z_{11})={{\rm e}^{-a_1-a_2}z_{12}}$
by using~\eqref{eq:cent}.
\end{proof}

\begin{proof}[Proof of Proposition~\ref{prop:key}]
Note that $\omega^{(\eta)}$ corresponds to $\omega^{-1}$
under the identification~\eqref{eq:eta-b}.
By using Proposition~\ref{prop:Gtheta}, \eqref{eq:Phi_1st_factor-2}, and
Proposition~\ref{prop:Phi_n_2nd_factor}, we have
\begin{align*}
\Phi_n
\bigl(\Gro{s_{\theta}}(z|\eta)\bigr)
&=
\Phi_n
\bigl(\Gro{s_1s_2\dots s_{n-1}}(z|\eta)\cdot \Gro{s_{n-2}\dots s_2s_1}\bigl(x|\omega^{(\eta)}(\eta)\bigr)
\bigr)\\
&=
\frac{\tau_n}{\tau_{n-1}}
\frac{{\rm e}^{a_n}}{z_{nn}}
\cdot
\omega^{-1}
\left(
\frac{{\rm e}^{-a_1-a_2}z_{12}}{\tau_1}
\right)\\
&=
\frac{\tau_n}{\tau_{n-1}}
\frac{{\rm e}^{a_n}}{z_{nn}}
\cdot
\frac{{\rm e}^{-a_n-a_1}z_{n,n+1}}{\tau_1}
\qquad\text{
since $\omega(z_{ij})=z_{i+1,j+1}$}
\\
&=
\frac{\tau_n}{\tau_{n-1}\tau_1}
\frac{{\rm e}^{-a_1}z_{n,n+1}}{z_{nn}}.
\end{align*}
From~\eqref{eq:zij_in_KkSchur},
we have
$z_{n,n+1}=\omega^{-1}(z_{12})={\rm e}^{a_n}z_{nn}g_{s_0}^{(k)}(y|b)$.
Thus, together with~\eqref{eq:sigma(tg)}, we have
\[
\frac{{\rm e}^{-a_1}z_{n,n+1}}{z_{nn}}
={\rm e}^{-a_1+a_n}g_{s_0}^{(k)}(y|b)
=\sigma^{-1}\bigl(\tilde{g}_{s_0}^{(k)}(y|b)\bigr).
\]
In view of
$
\tilde{\Phi}_n\bigl(Q^{\theta^\vee}\bigr)=
\frac{\tau_n}{\tau_{n-1}\tau_1}$,
the proof is complete.
\end{proof}

Now we can prove a~crucial property of $\tilde{\Phi}_n$.
\begin{Proposition}\label{prop:Phi_D0} We have
$\tilde{\Phi}_n\circ D_0^{Q}=D_0\circ\tilde{\Phi}_n
$.
\end{Proposition}
\begin{proof}
Since $s_\theta^*=s_\theta$, we have $\tilde{\Phi}_n\circ T_\theta^Q
=T_\theta \circ \tilde{\Phi}_n$
and $\tilde{\Phi}_n\circ s_\theta
=s_\theta \circ \tilde{\Phi}_n$.
From Proposition~\ref{prop:key}, we have
\begin{align*}
\tilde{\Phi}_n\circ D_0^{Q}&=\tilde{\Phi}_n
\circ \bigl(T_\theta^{Q}+
Q^{-\theta^\vee}
\Gro{s_\theta}
s_\theta
\bigr)=\tilde{\Phi}_n\circ T_\theta^{Q}+
\tilde{\Phi}_n\circ\bigl(Q^{-\theta^\vee}
\Gro{s_\theta}(x,y)
s_\theta\bigr)\\
&=T_\theta\circ\tilde{\Phi}_n
+
\tilde{\Phi}_n\bigl(Q^{-\theta^\vee}
\Gro{s_\theta}(x,y)\bigr)
\tilde{\Phi}_n\circ
s_\theta=T_\theta\circ\tilde{\Phi}_n+
\tilde{g}_{s_0}^{(k)}(y|b)\cdot
\tilde{\Phi}_n\circ
s_\theta
\\
&=T_\theta\circ\tilde{\Phi}_n+
D_0(1)\cdot
s_\theta\circ \tilde{\Phi}_n=D_0\circ\tilde{\Phi}_n,
\end{align*}
where the last equality follows from the simply identity
$
D_0=T_\theta
+ D_0(1)
s_\theta$.
\end{proof}

\begin{Corollary}\label{cor:Phi_commutes_w_Dem}
Let \smash{$x\in \hat{W}_G^0$}. Then
$
\tilde{\Phi}_n\circ D_x^Q=
D_{x^{\omega_k}}
\circ \tilde{\Phi}_n$.
\end{Corollary}
\begin{proof}
This is an immediate consequence of Proposition~\ref{prop:Phi_D0} and~\eqref{eq:D_and_Phi}.
\end{proof}

\begin{proof}[Proof of Theorem~\ref{thm:main}]
 \begin{align*}
 \tilde{\Phi}_n\bigl(Q^\xi \Gro{w}\bigr)
&=
\tilde{\Phi}_n\bigl(D_x^Q(1)\bigr) \qquad \text{by Proposition~\ref{prop:Quantum_D(1)}}\\
&={D}_{x^{\omega_k}}\bigl(\tilde{\Phi}_n(1)\bigr) \qquad \text{by Corollary~\ref{cor:Phi_commutes_w_Dem}}\\
&={D}_{x^{\omega_k}}(1)=\tilde{g}_{x^{\omega_k}}^{(k)}(y|b).\tag*{\qed}
\end{align*} \renewcommand{\qed}{}
\end{proof}

\section[Determinantal formulas for the k-small g\_la\^{}(k)(y|b)]{Determinantal formulas for the $\boldsymbol{k}$-small $\boldsymbol{\tilde{g}_\la^{(k)}(y|b)}$}\label{sec:det}

\subsection[S\_n-and T\_i-actions on O(Z)]{$\boldsymbol{S_n}$-and $\boldsymbol{T_i}$-actions on $\boldsymbol{\mathscr{O}(\mathscr{Z})}$}

Let $s_i$ $(1\le i\le n-1)$ be the generators of the symmetric group $S_n$, which acts on $\mathscr{O}\bigl(\tilde{\mathscr{Z}}\bigr)^{\reg}$ as a~$\CC$-algebra morphism by
\[
s_i\bigl({\rm e}^{\pm a_j}\bigr)={\rm e}^{\pm a_{s_i(j)}},\qquad
s_i(Z)=\mathfrak{s}_iZ\mathfrak{s}_i^{-1},
\qquad \mbox{where} \quad \mathfrak{s}_i=E-\bigl({\rm e}^{-a_i}-{\rm e}^{-a_{i+1}}\bigr)E_{i+1,i}.
\]

The $S_n$-action preserves the subalgebra $\mathscr{O}({\mathscr{Z}})$.
The isomorphism \smash{$\beta\colon\mathscr{O}(\mathscr{Z})\to \hat{\Lambda}_{(n)}^\Rep$} given in Theorem~\ref{thm:ISY} satisfies $s_i\circ\beta=\beta\circ s_i$ for all $i$ because $s_i(Z)$ is contained in the centralizer of~${s_i(A):=A|_{a_i\leftrightarrow a_{i+1}}}$.

On the coordinate functions $z_{ij}$, $S_n$ acts as
\begin{gather}
s_i(z_{jj})=z_{s_i(j),s_i(j)}\label{eq:s-2},\\
s_i(z_{ki})=z_{ki}+\bigl({\rm e}^{-a_i}-{\rm e}^{-a_{i+1}}\bigr)z_{k,i+1},\qquad k<i,\label{eq:s-3}\\
s_i(z_{i+1,j})=z_{i+1,j}-\bigl({\rm e}^{-a_i}-{\rm e}^{-a_{i+1}}\bigr)z_{ij},\qquad j>i,\label{eq:s-4}\\
s_i(z_{kj})=z_{kj},\qquad
k\ne i+1 \quad\text{or}\quad j\ne i.
\label{eq:s-5}
\end{gather}
It should be noted that~\eqref{eq:s-2} is a~consequence of~\eqref{eq:s-3}, and $z_{ij}$ is defined by~\eqref{eq:cent}.

\begin{Proposition}\label{prop:S_n_invariant_functions}
The action variables $c_0,c_1,\dots,c_{n-1}$ are $S_n$-invariant.
\end{Proposition}
\begin{proof}
From~\eqref{eq:def_of_c_i} and the definition of $s_i$, we have
\begin{align*}
s_i(Z)&=\mathfrak{s}_iZ\mathfrak{s}_i^{-1}
=\mathfrak{s}_i\left(\sum_{j=0}^{n-1}c_jA^j\right)\mathfrak{s}_i^{-1}\\
&=\sum_{j=0}^{n-1}c_j\mathfrak{s}_iA^j\mathfrak{s}_i^{-1}=\sum_{j=0}^{n-1}c_j\bigl(\mathfrak{s}_iA\mathfrak{s}_i^{-1}\bigr)^j=\sum_{j=0}^{n-1}c_js_i(A)^j.
\end{align*}
Since $s_i(A)$ is also regular, this equation uniquely determines $c_j$'s.
On the other hand, because~$s_i$ is a~$\C$-algebra homomorphism we have
$s_i(Z)=\sum_{j=0}^{n-1}s_i(c_j)s_i(A)^j$.
Hence we have $s_i(c_j)=c_j$.
\end{proof}
Theorem~\ref{thm:det_c_to_tausigma} implies the following important consequence.
\begin{Corollary}\label{cor:tausigma invariant}
 $\tau_i$ and $\sigma_i$ are $S_n$-invariant.
\end{Corollary}
\begin{Remark}\label{rem:PinvZP}
We have
\[
P^{-1}ZP=\sum_{j=0}^{n-1}c_j C_A^j.
\]
Since the entries of $C_A$ are $S_n$-invariant,
the entries of $P^{-1}ZP$ are $S_n$-invariant by Proposition~\ref{prop:S_n_invariant_functions}.
\end{Remark}

\begin{Proposition}\label{prop:T_action_on_z}
For $1\le i\le n-1$,
\begin{gather}
T_i(z_{ki})={\rm e}^{-a_i}z_{k,i+1},\label{eq:z_Tz}\\
T_i(z_{i+1,j})=-{\rm e}^{-a_i}z_{ij},\label{eq:z_Tz_2}\\
T_i(z_{kj})=0,\qquad
k\ne i+1 \quad \text{or}\quad j\ne i.\label{eq:z_Tz_3}
\end{gather}
\end{Proposition}

\begin{proof}
Equations~\eqref{eq:z_Tz}, \eqref{eq:z_Tz_2}, and~\eqref{eq:z_Tz_3} are immediate consequences of~\eqref{eq:s-3}, \eqref{eq:s-4}, and~\eqref{eq:s-5}, respectively.
\end{proof}

For convenience, we extend the definition of $z_{ij}$ to any $1\le i\le j$ by letting $z_{ii}=z_{i+n,i+n}=z_{i+2n,i+2n}=\cdots$ and
\[
z_{ij}=-
\frac{z_{i,j-1}-z_{i+1,j}}{{\rm e}^{-a_{i(\mathrm{mod}\, n)}}-{\rm e}^{-a_{j(\mathrm{mod}\, n)}}}.
\]
By definition, $z_{ij}$ is an element of $\mathscr{O}\bigl(\Tilde{\mathscr{Z}}\bigr)^{\reg}$.
Let $\omega:=s_1s_2\dots s_{n-1} $.
It is shown by induction on~${j-i\geq 0}$ that
\begin{equation}\label{eq:pi_zij}
\omega(z_{ij})=z_{i+1,j+1},
\end{equation}
which implies $z_{ij}\in \mathscr{O}(\mathscr{Z})$.
In particular, we have $z_{ij}=z_{i+n,j+n}$.

\begin{Remark}
$T_\theta$ can be seen as the $n$-th divided difference operator.
Proposition~\ref{prop:T_action_on_z} is naturally extended to
the case for $i=n$ as
$
T_\theta(z_{k,n})={\rm e}^{-a_n}z_{k,n+1}$ and $ T_\theta(z_{n+1,j})=-{\rm e}^{-a_n}z_{nj}$.
Moreover, the expressions~\eqref{eq:s-2}--\eqref{eq:z_Tz_2} are also valid for arbitrarily $i,k,j\in \ZZ$ under the identification $T_{i+n}=T_i$ and $T_n=T_\theta$.
\end{Remark}

\subsection[Determinantal formula for g\^{}(k)\_lambda]{Determinantal formula for $\boldsymbol{\tilde{g}^{(k)}_\lambda}$}
\label{sec:det_formula}
A $k$-bounded partition $\lambda$ is \textit{$k$-small} if $\lambda$ is contained in at least one of $R_1,R_2,\dots,R_{n-1}$.
This is equivalent to $\ell(\la)+\la_1\le n$.

\subsubsection{Notation}
We identify the set $\Iaf=\{0,1,\dots,n-1\}$ of
the type \smash{$A_{n-1}^{(1)}$} affine Dynkin nodes
with $\Z/n\Z$.
The $n$-\emph{residue} is the map
$
\mathrm{res}\colon
\mathbb{N}\times \mathbb{N}
\rightarrow \Iaf\cong \Z/n\Z$, $
(i,j)\to j-i \mod n$.
For $x=(i,j)\in \la$, let~${
\mathfrak{d}(x):=\{{\rm e}^{-a_{\mathrm{res}(s,j)}}\mid i\le s\le \la_j'\}}$.

For any subset $X$ of $\{{\rm e}^{-a_1},\dots,{\rm e}^{-a_n}\}$
we denote by $h_m(X)$ the $m$-th complete symmetric
polynomial in $X$.
We use
the abbreviation
\smash{$
\ha_m^{\langle i,j\rangle;\la}:=h_m(\mathfrak{d}(i,j))$},
for $(i,j)\in \la$.
For $k$-small partition $\la$, define
\begin{equation}
\label{eq:def xi}
\co_\la(y)=\prod_{i=1}^{\la_1'}\Omega(b_{\mathrm{res}(i,1)}|y).
\end{equation}

\begin{Example}\label{ex:h d}
When $n=6$ and
$
\lambda=(3,3,1)$.
If we fill the boxes of $\la$ with the $n$-residue, we have
\[
\begin{ytableau}
0&1&2\\
5&0&1\\4
\end{ytableau}.
\]
For example, $\mathfrak{d}(2,1)=\big\{{\rm e}^{-a_5},{\rm e}^{-a_4}\big\}$ so
\smash{$\ha_m^{\langle 2,1\rangle;\la}=h_m({\rm e}^{-a_5},{\rm e}^{-a_4})$}.
We also have
$
\co_\la(y)=\Omega(b_6|y)\allowbreak\times\Omega(b_5|y)\Omega(b_4|y)$.
\end{Example}

\subsubsection[Determinantal formula for g\_la(y|b) for k-small la]{Determinantal formula for $\boldsymbol{\tilde{g}_\la(y|b)}$ for $\boldsymbol{k}$-small $\boldsymbol{\la}$}

\begin{Theorem}\label{thm:determinant_formula_for_ksmall_g}
Let $\lambda$ be a~$k$-small partition.
Set $
l=n-\la_1'+1$
and $r=\la_1$.
Define the following
square matrix of size $n-l+1+r$
\begin{equation}\label{eq:M la}
 M_\la=
 \left(
\begin{array}{@{}cccc|cccc@{}}
 z_{ll}&\cdots &\cdots &z_{ln}&z_{l,n+1}&\cdots&\cdots&z_{l,n+r}\\
 &\ddots& &\vdots&\vdots&\cdots&\cdots&\vdots\\
 &&\ddots &\vdots&\vdots&\cdots&\cdots&\vdots\\
 & &&z_{nn}&z_{n,n+1}&\cdots&\cdots&z_{n,n+r}\\ \hline
 &\cdots
 & \ha_2^{\langle {2,1}\rangle;\la} & -\ha_1^{\langle 1,1\rangle;\la} & 1 & & \\
 &&\cdots
 &\ha_2^{\langle 2,2\rangle;\la}&-\ha_1^{\langle 1,2\rangle;\la} & \ddots &&\\
 &&\cdots&\vdots&\vdots & \ddots &\ddots&\\
 &&&\pm \ha_r^{\langle r,r\rangle;\la}&\mp \ha_{r-1}^{\langle r-1,r\rangle;\la}&\cdots &-\ha_1^{\langle 1,r\rangle;\la} &1
\end{array}
\right).
\end{equation}
Then we have
\begin{equation}\label{eq:matrix_M}
\tilde{g}^{(k)}_{\lambda}(y|b)=
\frac{\det(M_\lambda)}{\co_\la(y)}.
\end{equation}
\end{Theorem}

\begin{Remark}\quad
\begin{enumerate}
\itemsep=0pt
\item[(1)] $l$ is the $n$-residue of the bottom box of the first column of $\la$.
 \item[(2)] We have $\la_1<l$ because $\la$ is $k$-small.
If we take $r$ such that $\la_1\le r<l$ and
 consider the matrix $M_\la$ by the same formula~\eqref{eq:M la},
 then~\eqref{eq:matrix_M} also holds.
\end{enumerate}

\end{Remark}

\begin{Example}\label{ex:k small det}
Let $n=6$ and $\la$ be as in Example~\ref{ex:h d}. We have $l=4$.
The formula for $\tilde{g}^{(5)}_{\lambda}(y|b)$ by taking $r=3$ is
\[
\frac{1}{\Omega(b_4|y)\Omega(b_5|y)\Omega(b_6|y)}
\begin{vmatrix}
 z_{44}&z_{45}&z_{46}&z_{47}&z_{48}&z_{49}\\
 0&z_{55}&z_{56}&z_{57}&z_{58}&z_{59}\\
 0&0&z_{66}&z_{67}&z_{68}&z_{69}\\
 -f_3^{\langle 3,1\rangle;\la} & f_2^{\langle 2,1\rangle;\la} & -f_1^{\langle 1,1 \rangle;\la} & 1 & 0 & 0\\
 0&0&f_2^{\langle 2,2\rangle;\la}&-f_1^{\langle 1,2 \rangle;\la} & 1 &0\\
 0&0&0&f_2^{\langle 2,3\rangle;\la}&-f_1^{\langle 1,3 \rangle;\la} & 1
\end{vmatrix}.
\]
\end{Example}

\begin{Remark}
\bgroup
\renewcommand\binom[2]{\genfrac{[}{]}{0pt}{1}{#1}{#2}}
In the non-equivariant case $a_i=0$, \eqref{eq:matrix_M} reduces to
\[
\Tilde{g}^{(k)}_\lambda(y|0)=
\begin{vmatrix}
h_0(y)&h_1(y)&h_2(y)&\cdots&\cdots&\cdots &\cdots &h_{n-1}(y)\\
&h_0(y)&h_1(y)&h_2(y)&\cdots&\cdots &\cdots &h_{n-2}(y)\\
&&\ddots&\ddots&\ddots&\cdots &\cdots &\vdots\\
&&&h_0(y)&h_1(y)&h_2(y)&\dots &h_r(y)\\
\binom{\lambda_1'}{\lambda_1'}&\cdots&\binom{\lambda_1'}{2}&\binom{\lambda_1'}{1}&\binom{\lambda_1'}{0}\\
&\binom{\lambda_2'}{\lambda_2'}&\cdots&\binom{\lambda_2'}{2}&\binom{\lambda'_2}{1}&\binom{\lambda_2'}{0}\\
&&\ddots&\cdots&\cdots&\ddots&\ddots\\
&&&&\binom{\lambda_r'}{\lambda_r'}&\cdots&\binom{\lambda_r'}{1}&\binom{\lambda_r'}{0}
\end{vmatrix},
\]
where $h_i=h_i(y)$ and $\binom{a}{b}=(-1)^b\genfrac{(}{)}{0pt}{1}{a}{b}$.
\egroup
\end{Remark}

The proof of Theorem~\ref{thm:determinant_formula_for_ksmall_g} is divided into two cases: (i) the case when $\lambda$ is a~one-column partition and (ii) the case when $\lambda$ is a~general $k$-small partition.

\subsubsection[Case: la=(1\^{}i)]{Case: $\boldsymbol{\la=\bigl(1^i\bigr)}$}
Let $\lambda$ is a~one-column partition $\lambda=\bigl(1^i\bigr)$.
We have
\[
 M_{(1^i)}=
 \left(
\begin{array}{@{}cccc|c@{}}
 z_{n-i+1,n-i+1}&\cdots &\cdots &z_{n-i+1,n}&z_{n-i+1,n+1}\\
 &\ddots& &\vdots&\vdots\\
 &&\ddots &\vdots&\vdots\\
 & &&z_{nn}&z_{n,n+1}\\ \hline
 (-1)^i\ha_i^{\langle {i,1}\rangle;(1^i)} &\cdots
 & \ha_2^{\langle {2,1}\rangle;(1^i)} & -\ha_1^{\langle 1,1\rangle;(1^i)}&1
\end{array}
\right).
\]

In this case, \smash{$\det(M_{(1^i)})$} given in~\eqref{eq:matrix_M} is expanded as
\[
\det(M_{(1^i)})
=
\sum_{m=0}^{i}
\ha_{m}^{\langle m,1\rangle;(1^i)}
d^{(i)}_m \\
\]
by the expansion along the bottom row,
where
\smash{$
d^{(i)}_m=\det Z_{[n-i+1,n]}
^{[n-i+1,n+1]\setminus \{n-m+1\}}$}.
For $0\le m\le i< n$, define
\begin{gather}
f^{(i)}_m=
\ha_{m}^{\langle m,1\rangle; (1^i)}
\cdot d^{(i)}_m.\label{eq: varphi m}
\end{gather}
Then~\eqref{eq:matrix_M} reads
\begin{equation}\label{eq:$K$-$k$-one_column_determinant}
\begin{aligned}
\tilde{g}^{(k)}_{(1^i)}(y|b)
=
\frac{
\sum_{m=0}^{i}
\varphi_m^{(i)}
}{\co_{(1^i)}(y)}
\end{aligned}.
\end{equation}

We consider the minor
$d^{(i+1)}_m$
of size $i+1$.
Note that we have
\begin{equation}
 d^{(i+1)}_m=z_{n-i,n-i}d^{(i)}_{m} \label{eq: d=dz}
\end{equation}
since $Z$ is upper triangular.
It follows that
\begin{equation}\label{eq:tilde varphi}
z_{n-i,n-i}\varphi_m^{(i)}=\ha_{m}^{\langle m,1\rangle; (1^i)}d^{(i+1)}_m.
\end{equation}

\begin{Lemma}\label{lem:T_on_f}
For $1\le i\le n-1$, we have \smash{$T_{n-i}\bigl(d^{(i+1)}_i\bigr)={\rm e}^{-a_{n-i}}d^{(i+1)}_{i+1}$}.
\end{Lemma}
\begin{proof}
From~\eqref{eq:z_Tz_2} and~\eqref{eq:z_Tz_3}, we have
\begin{equation}
\label{eq:T d i}
T_{n-i}\bigl(d^{(i)}_{i}\bigr)
=
T_{n-i}\left(\det Z_{[n-i+1,n]}^{[n-i+2,n+1]}\right)=-{\rm e}^{-a_{n-i}}\cdot \det Z_{[n-i,n]\setminus \{n-i+1\}}^{[n-i+2,n+1]}
\end{equation}
because every entry of the rows of \smash{$Z_{[n-i+1,n]}^{[n-i+2,n+1]}$} except for the first row is invariant under $s_{n-i}$.
Then by using this, we have
\begin{align*}
T_{n-i}\bigl(d^{(i+1)}_{i}\bigr)={}&T_{n-i}\bigl(z_{n-i,n-i}d^{(i)}_{i}\bigr)\qquad\text{by~\eqref{eq: d=dz}}\\
={}&T_{n-i}(z_{n-i,n-i})\cdot d^{(i)}_{i}
+s_{n-i}(z_{n-i,n-i})\cdot T_{n-i}\bigl(d^{(i)}_{i}\bigr) \\
={}&
{\rm e}^{-a_{n-i}}z_{n-i,n-i+1}\cdot
d^{(i)}_i
\\
&-z_{n-i+1,n-i+1}\cdot {\rm e}^{-a_{n-i}}
\det Z_{[n-i,n]\setminus \{n-i+1\}}^{[n-i+2,n+1]}\quad\text{by~\eqref{eq:z_Tz} and~\eqref{eq:T d i}}\\
={}&{\rm e}^{-a_{n-i}}d^{(i+1)}_{i+1},
\end{align*}
where the last equality follows by the expansion along the $1$st column of $d^{(i+1)}_{i+1}$.
\end{proof}
\begin{Lemma}\label{lem: D z varphi}
For $1\le i\le n-1$,
and $0\le m\le i$, we have
\[
D_{n-i}\bigl(z_{n-i,n-i}\varphi^{(i)}_m\bigr)=
\begin{cases}
\varphi^{(i+1)}_{m}, & m<i,\\
\varphi^{(i+1)}_{i}+\varphi^{(i+1)}_{i+1}, & m=i.
\end{cases}
\]
\end{Lemma}
\begin{proof}

Suppose $m<i$. We claim that $d^{(i+1)}_m$ is invariant under $s_{n-i}$.
In fact, we have
\begin{equation*}
 d^{(i+1)}_m=z_{n-i,n-i}z_{n-i+1,n-i+1}d^{(i-1)}_m
\end{equation*} by~\eqref{eq: d=dz},
which is $s_{n-1}$-invariant
by~\eqref{eq:s-2} and~\eqref{eq:z_Tz_3}.
We have
\begin{align*}
D_{n-i}\bigl({z_{n-i,n-i}{\varphi}}^{(i)}_m\bigr)&=
D_{n-i}\bigl(\ha_{m}^{\langle m,1\rangle; \bigl(1^i\bigr)}\cdot d^{(i+1)}_m\bigr)\qquad\text{by $\eqref{eq:tilde varphi}$}\\
&=D_{n-i}\bigl(\ha_{m}^{\langle m,1\rangle; \bigl(1^i\bigr)}\bigr)\cdot d^{(i+1)}_m\qquad\text{since $d^{(i+1)}_m$ is $s_{n-i}$-invariant}
\\&=
\ha_{m}^{\langle m,1\rangle; (1^{i+1})}\cdot d^{(i+1)}_m
\qquad\text{by Lemma~\ref{lem:T_on_h} }
\\&=\varphi^{(i+1)}_m.
\end{align*}

When $m=i$, noting that
$\ha_i^{\langle i,1\rangle; (1^{i})}
={\rm e}^{-ia_{n-i+1}}$,
we have
\begin{align*}
D_{n-i}\bigl({z_{n-i,n-i}{\varphi}}^{(i)}_i\bigr)
={}&
D_{n-i}\bigl(
\ha_i^{\langle i,1\rangle; (1^{i})}
d^{(i+1)}_i\bigr)\qquad\text{by~\eqref{eq:tilde varphi}}\\
={}&
D_{n-i}\bigl(\ha_i^{\langle i,1\rangle; (1^{i})}\bigr)\cdot d^{(i+1)}_i
+{\rm e}^{-ia_{n-i}}\cdot T_{n-i}\bigl(d^{(i+1)}_i\bigr)\\
={}&
h_i({\rm e}^{-a_{n-i}},{\rm e}^{-a_{n-i+1}})\cdot d_{i}^{(i+1)}\\
&+{\rm e}^{-(i+1)a_{n-i}}d^{(i+1)}_{i+1}
\qquad\text{by Lemmas~\ref{lem:T_on_h} and~\ref{lem:T_on_f}}\\
={}&\ha_i^{\langle i,1\rangle; (1^{i+1})}\cdot d_{i}^{(i+1)}+
\ha_{i+1}^{\langle i+1,1\rangle;(1^{i+1})}\cdot d^{(i+1)}_{i+1}
\\
={}&\varphi^{(i+1)}_i+\varphi^{(i+1)}_{i+1} \qquad\text{by~\eqref{eq: varphi m}}.\tag*{\qed}
\end{align*} \renewcommand{\qed}{}
\end{proof}

\begin{proof}[Proof of Theorem~\ref{thm:determinant_formula_for_ksmall_g} for $\boldsymbol{\la=(1)^i}$]
We show~\eqref{eq:$K$-$k$-one_column_determinant} by induction on $i\geq 1$.
When $i=1$, the desired equation is shown directly as follows:
\begin{align*}
\tilde{g}^{(k)}_{(1)}(y|b)
&=D_0(1)=
1+\Omega(b_1|y) T_\theta\left(\frac{1}{\Omega(b_1|y)}\right)
=
1+\Omega(b_1|y) T_\theta\left(\frac{z_{nn}}{\Omega(b_1|y)\Omega(b_n|y)}\right)\\
&=1+\frac{T_\theta(z_{nn})}{\Omega(b_n|y)}
=\frac{z_{nn}+{\rm e}^{-a_n}z_{n,n+1}}{\Omega(b_n|y)}
=\frac{\varphi^{(1)}_0+\varphi^{(1)}_1}{\Omega(b_n|y)}
=\frac{\varphi^{(1)}_0+\varphi^{(1)}_1}{\co_{(1)}(y)}.
\end{align*}
Suppose~\eqref{eq:$K$-$k$-one_column_determinant} holds for some $1\le i\le n-2$.
By induction hypothesis,
\begin{align*}
\tilde{g}^{(k)}_{(1^{i+1})}(y|b)
&=D_{n-i}\left(
\frac{\sum_{m=0}^{i}\varphi^{(i)}_m}{\co_{(1^i)}(y)}
\right)=
D_{n-i}\left(
\frac{z_{n-i,n-i}\sum_{m=0}^{i}{{\varphi}}^{(i)}_m}{\Omega(b_{n-i}|y)\co_{(1^i)}(y)} \right)
\\
&=
\frac{
\sum_{m=0}^{i}D_{n-i}\left({z_{n-i,n-i}{\varphi}}^{(i)}_m\right) }{\Omega(b_{n-i}|y)\co_{(1^i)}(y)}=
\frac{
\sum_{m=0}^{i+1}\varphi^{(i+1)}_m }{\co_{(1^{i+1})}(y)}\qquad\text{by Lemma~\ref{lem: D z varphi}},
\end{align*}
where in the third equality we use
fact that $\Omega(b_{n-i}|y)\co_{(1^i)}(y)$ is $s_{n-i}$ invariant.
\end{proof}

\subsubsection[Case: general k-small la]{Case: general $\boldsymbol{k}$-small $\boldsymbol{\la}$}
Before proving Theorem~\ref{thm:determinant_formula_for_ksmall_g} for general $k$-small $\lambda$, we list a~few important properties of ${M}_\lambda$.

\begin{Lemma}\label{lem:comparing_determinants}
Let $\lambda$ be a~$k$-small partition.
Assume $\la$ is $i$-addable.
Let $\kappa$ be a~partition
obtained from $\lambda$ by adding a~box of $n$-residue $i$. We also assume $\kappa$ is $k$-small.
Then
$
 D_i(\det M_\lambda)=\det M_\kappa $ if $i\neq 0$, and $
D_\theta(\det M_\lambda)=\det M_\kappa$.
\end{Lemma}

\begin{proof}
We may assume $\la\ne \varnothing$.
The integers $l$, $r$ are defined in Theorem~\ref{thm:determinant_formula_for_ksmall_g}.
Let $\hat{x}:=x-l+1$ and $x^\dagger:=x-l+1+n$ for $x\in \ZZ$.
In $M_\lambda$, $z_{q,p}$ is at the $(\hat{q},\hat{p})$-th position (if exists).
For~${1\le s\le r}$, the
$\bigl(s^\dagger,\hat{p}\bigr)$-th entry of $M_\lambda$ is either a~constant or a~complete symmetric polynomial in ${\rm e}^{-a_{b}},{\rm e}^{-a_{b+1}},\dots,{\rm e}^{-a_{p}}$ where $b$ is the $n$-residue of the bottom box of the $s$-th column.

Let \smash{$M_\lambda^{(i)}$} be the matrix obtained from $M_\lambda$ by applying the following column or row operations:
\begin{itemize}\itemsep=0pt
 \item[(a)] Add ${\rm e}^{-a_{i}}$ times the $\big(\hat{i}+1\big)$-th column to the $\hat{i}$-th column.
 \item[(b)] If $l\le i<n$, subtract ${\rm e}^{-a_{i}}$ times the $\hat{i}$-th row from the $\big(\hat{i}+1\big)$-th row.
\end{itemize}
Since $\lambda$ is $k$-small and $\la$ is $i$-addable, there is
a unique integer $j$ such that the $n$-residue of the bottom box of the $j$-th column is $i + 1$.
In \smash{$M^{(i)}_\lambda$}, all the entries in rows except for the $j^\dagger$-th one are $s_{i}$-invariant because:
\begin{itemize}
\itemsep=0pt
 \item the $(\hat{q},\hat{p})$-th entry is
 \begin{gather*}
 -{\rm e}^{a_i-a_{i+1}}z_{i,i+1}\qquad \text{if} \ q=i+1,\quad p=i,\\
 z_{qi}+{\rm e}^{-a_i}z_{q,i+1}\qquad \text{if} \ q\neq i+1,\quad p=i,\\
 z_{i+1,p}-{\rm e}^{-a_i}z_{ip} \qquad \text{if} \ q=i+1,\quad p\neq i,\\
 z_{qp} \qquad \text{if} \ q\neq i+1,\quad p\neq i,
 \end{gather*}
One can verify that these entries are $s_i$-invariant from~\eqref{eq:s-2}--\eqref{eq:s-5}.
\item the $\bigl(s^\dagger,\hat{i}\bigr)$-th entry $(s\neq j)$ is either a~constant or a~polynomial of the form
\[
h_t({\rm e}^{-a_{b}},\dots,{\rm e}^{-a_i})-{\rm e}^{-a_i}h_{t-1}({\rm e}^{-a_{b}},\dots,{\rm e}^{-a_i},{\rm e}^{-a_{i+1}}),
\]
up to a~sign,
where $b$ is the $n$-residue of the bottom box of the $s$-th column, and $t>0$.
It is not difficult to show the polynomial above is
$s_i$-invariant.

\item the $\bigl(s^\dagger,\hat{p}\bigr)$-th entry $(s\neq j,\ p\neq i)$ is either a~constant or a~symmetric polynomial in~\smash{${\rm e}^{-a_{b}},{\rm e}^{-a_{b+1}},\dots,{\rm e}^{-a_p}$} with $b$ the $n$-residue of the bottom box of
the $s$-th column. Note that we have $b\ne i+1$.
Clearly, these are $s_i$-invariant.
\end{itemize}
Therefore, it suffices to compute the $D_i$-actions on the $j^\dagger$-th row of \smash{$M^{(i)}_\lambda$} only.
Set
$m=\lambda'_j$.
Then
$m=j-i+n-1$,
and the \smash{$\bigl(j^\dagger,\hat{p}\bigr)$}-th entry of \smash{$M^{(i)}_\lambda$} is
\begin{gather*}
(-1)^{m-p+i+1}h_{m-p+i+1}\bigl({\rm e}^{-a_{i+1}},\dots,{\rm e}^{-a_p}\bigr) \qquad \text{if} \ \hat{p}<\hat{i},\\
(-1)^{m+1}{\rm e}^{-a_i}h_{m}\bigl({\rm e}^{-a_{i+1}}\bigr) \qquad \text{if} \ \hat{p}=\hat{i},\\
0 \qquad \text{if} \ \hat{p}>\hat{i}
\end{gather*}
and the $\bigl(j^\dagger,\hat{p}\bigr)$-th entry of $M^{(i)}_\kappa$ is
\begin{gather*}
(-1)^{m-p+i+1}h_{m-p+i+1}({\rm e}^{-a_i},{\rm e}^{-a_{i+1}},\dots,{\rm e}^{-a_p}) \qquad \text{if} \
\hat{p}<\hat{i},\\
(-1)^{m+1}
(
h_{m+1}({\rm e}^{-a_i})
-{\rm e}^{-a_i}h_{m}({\rm e}^{-a_i},{\rm e}^{-a_{i+1}})
) \qquad \text{if} \  \hat{p}=\hat{i},\\
0   \qquad \text{if} \ \hat{p}>\hat{i}.
\end{gather*}
By Lemma~\ref{lem:T_on_h}, the image of the $j^\dagger$-th row of \smash{$M_\lambda^{(i)}$} under $D_i$ coincides with the $j^\dagger$-th row of~\smash{$M_\kappa^{(i)}$}.
\end{proof}

\begin{proof}[Proof of Theorem~\ref{thm:determinant_formula_for_ksmall_g}]
Let $w$ be the element of \smash{$\hat{W}_{G}^0$} associated with $\lambda$.
Since $\lambda$ is $k$-small, there exists a~reduced expression $w=vs_ls_{l+1}\dots s_{n-1}s_0$ for some \label{$w'$ is not affine Grassmannian.}
$v\in \hat{W}_{G}$, where $l$ is the $n$-residue at the bottom of the first column of $\lambda$.

We show the theorem by induction on $\ell(v)$.
The case $\ell(v)=0$ is done, so suppose $\ell(v)>0$.
In this case, $\lambda$ contains at least one removable box
not included in the first column.
Let $i$ be the $n$-residue of the box.
Let $\mu$ the partition obtained from $\la$ by removing
the box of $n$-residue $i$.

(i) If $i\neq 0$, by induction hypothesis, we have
\[
\tilde{g}^{(k)}_\lambda(y|b)=
D_i\bigl(
\tilde{g}^{(k)}_{\mu}(y|b)
\bigr)
=
D_i
\left(
\frac{
\det(M_{\mu})
}{\co_\la(y)}
\right).
\]
Since $\lambda$ is $k$-small, we easily show that $\co_\la(y)$ is invariant under $s_i$.
Therefore, from Lemma~\ref{lem:comparing_determinants}, we have
\[
D_{i}\left(
\frac{\det M_{\mu}}{\co_\la(y)}
\right)
=
\frac{
D_{i}(\det M_{\mu})}{\co_\la(y)}
=
\frac{\det M_{\lambda}}{\co_\la(y)}.
\]

(ii) For $i=0$, because $\Omega(b_1|y)\co_\la(y)$ is $s_\theta$-invariant, we have
\[
\begin{aligned}
\tilde{g}^{(k)}_\lambda(y|b)
&
=
\Omega(b_1|y)
D_\theta
\left(
\frac{
\det M_{\mu}
}{\Omega(b_1|y)\co_\la(y)}
\right)\\
&
=
\Omega(b_1|y)
\frac{
D_\theta(\det M_{\mu})}{\Omega(b_1|y)\co_\la(y)}
=
\frac{\det M_{\lambda}}{\co_\la(y)}
\end{aligned}
\]
from Lemma~\ref{lem:comparing_determinants}.
\end{proof}

\begin{Corollary}\label{cor:g_rectangle}
Suppose $i+j\le n$.
Then we have
\[
\tilde{g}^{(k)}_{(i^j)}(y|b)=
\frac{
\det(M'_{(i^j)})
}{\co_{(i^j)}(y)},\qquad
M'_{(i^j)}=\omega^{-j}
\left(
\begin{matrix}
Z_{[1,j]}^{[1,i+j]}\\\hline
{
\bigl(P^{-1}\bigr)}^{[1,i+j]}_{[j+1,i+j]}
\end{matrix}
\right).
\]
\end{Corollary}

\begin{proof}
Note also that the size of the square matrix $M_{(i^j)}$ is $i+j$.
Let $W$ be the matrix consisting of the last $i$ rows of $M_{(i^j)}$.
Then we can write \begin{equation}
M_{(i^j)}=\left(
\begin{matrix}
\omega^{-j}Z_{[1,j]}^{[1,i+j]}\\\hline
W
\end{matrix}
\right).\label{eq:M rec}
\end{equation}

There is a~lower unitriangular $i\times i$ matrix $N$
such that
\begin{equation}
N \omega^j(W)
=\bigl(P^{-1}\bigr)_{[j+1,i+j]}^{[1,i+j]}.
\label{eq:N omega j W}
\end{equation}
In fact, if we define
$
N:
=
N_1
N_2
\cdots
N_{i-1}$,
where
\smash{$N_m=1-\sum_{j=1}^m{\rm e}^{-a_j}E_{r+j-m,r+j-m-1}$}, then by using the explicit formula for $P^{-1}$, it is straightforward to check~\eqref{eq:N omega j W}.

Combined with~\eqref{eq:M rec} and~\eqref{eq:N omega j W}, we have
\[
\begin{pmatrix}
E_j & O\\
O & N
\end{pmatrix}\omega^j M_{(i^j)}=\left(
\begin{matrix}
Z_{[1,j]}^{[1,i+j]}\\\hline
N\omega^j W
\end{matrix}
\right)=\left(
\begin{matrix}
Z_{[1,j]}^{[1,i+j]}\\\hline
\bigl(P^{-1}\bigr)_{[j+1,i+j]}^{[1,i+j]}
\end{matrix}
\right)=\omega^j M_{(i^j)}'.
\]
Therefore, by taking determinants we
obtain the desired result.
\end{proof}

\begin{Corollary}\label{cor:k rect and tau}
We have
 \begin{align}
 \tilde{g}^{(k)}_{R_i}(y|b)&=
 \frac{\sigma_{n-i}}{\co_{R_i}(y)},\label{eq:det sigma i}
 \\ g^{(k)}_{R_i}(y|b)&=
{\rm e}^{\sum_{s=0}^{n-i-1}a_{i-s}}
 \frac{
 \tau_{n-i}}{\co_{R_i}(y)}.\label{eq:det tau i}\end{align}
\end{Corollary}
\begin{proof}
Let $M'_{R_i}$ be the $n\times n$ matrix in Corollary~\ref{cor:g_rectangle}. We have
\[
\omega^{i}M'_{R_i}=
\left(
\begin{matrix}
Z_{[1,n-i]}^{[1,n]}\\\hline
{
\bigl(P^{-1}\bigr)}^{[1,n]}_{[n-i+1,n]}
\end{matrix}
\right).
\]
So we have
\[
\bigl(\omega^{i}M'_{R_i}\bigr)P=
\left(
\begin{matrix}
(ZP)_{[1,n-i]}^{[1,n]}\\
\hline
O_{n-i}|E_{i}
\end{matrix}
\right).
\]
By taking determinants, we have
\[
\omega^{i}\det(M'_{R_i})
=\det\bigl(\omega^{i}M'_{R_i}\bigr)=\det\left((ZP)^{[1,n-i]}_{[1,n-i]}\right)=\sigma_{n-i}.
\]
Since $\sigma_{n-i}$ is $S_n$-invariant (see Corollary~\ref{cor:tausigma invariant}), we obtain
\eqref{eq:det sigma i}.
It is straightforward to obtain~\eqref{eq:det tau i} from~\eqref{eq:det sigma i}
by using~\eqref{eq:sigma(Omega)} and
Proposition~\ref{prop:v1_of_kKShur}.
\end{proof}

\subsection[Determinantal formula for g\^{}(k)\_lambda]{Determinantal formula for $\boldsymbol{g^{(k)}_\lambda}$}

The determinantal formula for the $K$-$k$-Schur function $g^{(k)}_\lambda(y|b)$ is also obtained for a~$k$-small $\lambda$.
\begin{Theorem}\label{thm:determinantal_formula_for_g}
Let $\lambda$ be a~$k$-small partition.
Let $l$, $r$ be as in Theorem~{\rm\ref{thm:determinant_formula_for_ksmall_g}}.
Let $W_\la$ be the matrix consisting of the last $i$ rows of $M_{\la}$.
Then \begin{equation}
g^{(k)}_{\lambda}(y|b)
={\rm e}^{\sum_{x\in \mathrm{diag}(\lambda)}(a_{r(x)+1}-a_{b(x)})}
\frac{\det(N_\lambda)}{\co_\la(y)},\qquad
N_\lambda=
\left(
\begin{matrix}
 (ZA)_{[l,n]}^{[l,n+r]} \\ \hline
 W_{\la}
\end{matrix}
\right).
\label{eq:det g small la}
\end{equation}
\end{Theorem}
\begin{proof} Since $\sigma^{-1}(Z)=ZA$ and $\sigma$ is $\Rep$-linear, we have
$\sigma^{-1}(M_\la)=N_\la.
$ In view of this,
\eqref{eq:det g small la}
is obtained from~\eqref{eq:matrix_M} by applying~\eqref{eq:sigma(Omega)} and Proposition~\ref{prop:v1_of_kKShur}.
\end{proof}

\begin{Example}\label{eq:n6_7boxes} For $\la$ in Example
 \ref{ex:k small det}, $g^{(5)}_{\lambda}(y|b)$ is
\[
\begin{aligned}
&
\dfrac{{\rm e}^{a_2+a_3+a_5}}
{\Omega(b_4|y)\Omega(b_5|y)\Omega(b_6|y)}
\begin{vmatrix}
 z'_{44}&z'_{45}&z'_{46}&z'_{47}&z'_{48}&z'_{49}\\
 0&z'_{55}&z'_{56}&z'_{57}&z'_{58}&z'_{59}\\
 0&0&z'_{66}&z'_{67}&z'_{68}&z'_{69}\\
 -f_3^{\langle 3,1\rangle;\la} & f_2^{\langle 2,1\rangle;\la} & -f_1^{\langle 1,1 \rangle;\la} & 1 & 0 & 0\\
 0&0&f_2^{\langle 2,2\rangle;\la}&-f_1^{\langle 1,2 \rangle;\la} & 1 &0
\end{vmatrix},
\end{aligned}
\]
where
$z'_{ij}={\rm e}^{-a_i}z_{ij}-z_{i,j-1}$.
Note that
$z_{ij}'
$ is specialized to
$h_{j-i}(y)-h_{j-i-1}(y)$
at $a_i=0$.
\end{Example}

\section[k-rectangle factorization property]{$\boldsymbol{k}$-rectangle factorization property}\label{sec:krect}
Define for $1\le j\le n$
\begin{equation}
\label{eq:def de}
\de_j:=
\Omega(b_{n-j+1}|y)^{-1}\cdots \Omega(b_{n}|y)^{-1}.
\end{equation}
Note that for $1\le i\le n-1$, we have
\begin{equation}
\de_{n-i}
=\co_{R_{i}}(y)^{-1}.
\label{eq:d Omega R}
\end{equation}
\begin{Lemma}\label{lem:commutation_relation_Omega}

Let $1\le j\le n$. We consider $\de_j$ a~linear transformation of
$\hat{\Lambda}^{R(T)}$ by multiplication.
For any $i\in \tilde{I}$, we have
\begin{equation}
{D}_i \circ \de_{j}
=\de_j\circ \bigl(\omega^{-j}D_{i+j} \omega^{j}\bigr),
\label{eq:D d}
\end{equation}
where the subscripts are taken modulo $n$.
\end{Lemma}
\begin{proof}
Let us consider the case $i\ne 0$ and $i+j= 0\mod n $.
We have
\begin{align*}
D_i \circ \de_j&=
D_i \circ \Omega(b_{i+1}|y)^{-1}\dots \Omega(b_{n}|y)^{-1}\qquad
\text{by assumption $i+j\equiv 0 \mod n$}\\
&=D_i \circ
\Omega(b_{i}|y)
\Omega(b_{i}|y)^{-1}
\Omega(b_{i+1}|y)^{-1}\de_{j-1}\\
&=\Omega(b_{i}|y)^{-1}
\Omega(b_{i+1}|y)^{-1}\de_{j-1}
\cdot
D_i \circ
\Omega(b_{i}|y)=\de_{j}
\circ
\Omega(b_{i}|y)^{-1}
\circ
D_i \circ
\Omega(b_{i}|y).
\end{align*}
In the third equality, we used that
$\Omega(b_{i}|y)^{-1}
\Omega(b_{i+1}|y)^{-1}\de_{j-1}$
is $s_i$-invariant.
It is straightforward to show
$
\Omega(b_{i}|y)^{-1}
\circ
D_i \circ
\Omega(b_{i}|y)=\omega^i
\circ D_0 \circ \omega^{-i}$.
Hence we have~\eqref{eq:D d} in this case.
The other cases are left to the reader.
\end{proof}

\begin{Theorem}\label{thm:$k$-rect}
Let $\lambda$ be an arbitrary $k$-bounded partition.
Then we have
\[
\tilde{g}^{(k)}_{\lambda\cup R_i}(y|b)=
\tilde{g}^{(k)}_{\lambda}\bigl(y| \omega^{i}b\bigr)\cdot\tilde{g}^{(k)}_{R_i}(y|b).
\]
\end{Theorem}
\begin{proof}
Let $w_\lambda=s_{j_1}s_{j_2}\cdots s_{j_l}$ be a~reduced expression.
Let $w_\la'=s_{j_1+i}s_{j_2+i}\cdots s_{j_l+i}$.
It is shown in \cite[Lemma~2.15]{Tak} that we have a~length additive decomposition
$
w_{\la \cup R_i}=w'_\la w_{R_i}$.

Therefore, we have
\begin{align*}
\tilde{g}^{(k)}_{\lambda\cup R_i}(y|b)&=D_{w'_\la}D_{R_i}(1)=D_{w'_\la}\bigl(\tilde{g}^{(k)}_{R_i}(y|b)\bigr)
\\
&=D_{w'_\la}\left(\frac{\sigma_{n-i}}{\co_{R_i}(y)}\right)\qquad\text{by~\eqref{eq:det sigma i}}\\
&=\sigma_{n-i}\cdot D_{w'_\la}{\co_{R_i}(y)}^{-1}\qquad\text{by Proposition~\ref{prop:S_n_invariant_functions}}
\\
&=\sigma_{n-i}\cdot D_{w'_\la}\de_{n-i}\qquad\text{by~\eqref{eq:d Omega R}}\\
&=\sigma_{n-i} \de_{n-i}\cdot \omega^{i} D_{w_\la}\omega^{-i}(1)
\qquad\text{by Lemma~\ref{lem:commutation_relation_Omega}}
\\
&=\tilde{g}^{(k)}_{R_i}
(u|b)
\cdot \omega^{i}( D_{w_\la}(1))=\tilde{g}^{(k)}_{R_i}
(u|b)
\cdot \tilde{g}^{(k)}_{\la}
\bigl(u| \omega^i b\bigr).\tag*{\qed}
 \end{align*}\renewcommand{\qed}{}
\end{proof}

\section{Proof of Theorem~\ref{thm:max factor}}\label{sec:max factor}

This is obtained by computing the image of $\mathfrak{G}_{w_\circ}^Q$
under $\tilde{\Phi}_n$ in two different ways.

\subsection[Phi\_n(fac\_i)]{$\boldsymbol{\tilde{\Phi}_n(\fac_i)}$}

Recall that $\fac_i$ $(1\le i\le n-1)$ is defined by~\eqref{eq:def phi}.
\begin{Proposition}\label{prop:phi}
For $1\le i\le n-1$,
\[
\tilde{\Phi}_n(\fac_i)=
\frac{\prod_{l=1}^i\Omega(b_{i+l+1}|y)}{\sigma_i}\cdot\tilde{g}_{(n-i-1)^i}^{(k)}\bigl(y|\omega^{2i+1}b\bigr).
\]
\end{Proposition}
Note that for $i=n-1$,
the partition appearing on the right-hand side is $\varnothing$.
Since $\mathfrak{G}_{w_\circ}^Q$ is the product of $\fac_i$ $(1\le i\le n-1)$, we immediately obtain the following formula.
\begin{Corollary}\label{cor:image_of_Groth}
We have
\begin{equation}\label{eq:image_of_Groth}
\tilde{\Phi}_n\bigl(\mathfrak{G}_{w_\circ}^Q\bigr)
=
\frac{\prod_{i=1}^{n-1}\prod_{l=1}^i\Omega(b_{i+l+1}|y)}{{\sigma_1\sigma_2\cdots \sigma_{n-1}}}\cdot
{
\prod_{i=1}^{n-2}
\tilde{g}^{(k)}_{(n-i-1)^i}\bigl(y|\omega^{2i+1}b\bigr)
}.
\end{equation}
\end{Corollary}
Here we note
that the factor
$\prod_{i=1}^{n-1}\prod_{l=1}^i\Omega(b_{i+l+1}|y)$ in~\eqref{eq:image_of_Groth}
can be expressed as follows:
\begin{equation}
\label{eq:Omega_n}
\prod_{i=1}^n\Omega(b_i|y)^m\qquad \text{for} \ n=2m+1,\qquad
\prod_{i=1}^m\Omega(b_{2i-1}|y)^{m}\Omega(b_{2i}|y)^{m-1}\qquad \text{for} \  n=2m.
\end{equation}

As an application of Corollary~\ref{cor:image_of_Groth}, we have a~proof of Proposition~\ref{prop:long_iota_invariant}.
\begin{proof}[Proof of Proposition~\ref{prop:long_iota_invariant}]
The product \smash{$\prod_{i=1}^{n-2}\tilde{g}^{(k)}_{(n-i-1)^i}\bigl(y|\omega^{2i+1}b\bigr)$} is $\iota$-invariant because
\[
\iota\bigl(\tilde{g}^{(k)}_{(n-i-1)^i}\bigl(y|\omega^{2i+1}b\bigr)\bigr)
=
\tilde{g}^{(k)}_{i^{(n-i-1)}}\bigl(y|\omega^{-2i-1}b\bigr)
=\tilde{g}^{(k)}_{i^{(n-i-1)}}\bigl(y|\omega^{2(n-i-1)+1}b\bigr)
\]
 by
Proposition~\ref{prop:iota $K$-k}.
Let $\Omega_n$ denote the factor~\eqref{eq:Omega_n}.
We will show the factor $\Omega_n\cdot (\sigma_1\sigma_2\cdots \sigma_{n-1})^{-1}$ is also $\iota$-invariant.
In fact, from~\eqref{eq:Omega_n} and~\eqref{eq:iota Omega} we have
$
\iota(\Omega_n)
=\Omega_n/ \sigma_n^{n-1}.
$
From~\eqref{prop:S_on_tau},
we have~${
\iota(\sigma_1\cdots \sigma_{n-1})=\sigma_1\cdots \sigma_{n-1}\cdot \sigma_n^{n-1}}$.

Hence, from Corollary~\ref{cor:image_of_Groth},
$\tilde{\Phi}_n\bigl(\mathfrak{G}_{w_\circ}^Q\bigr)$ is invariant under the action of $\iota$.
\end{proof}

The rest of this subsection is devoted to the proof of Proposition~\ref{prop:phi}.

We start with the expression
\[
\tilde{\Phi}_n(\fac_i)=
\frac{\det\left((E-{\rm e}^{a_{i+1}}C_A )P^{-1}ZP\right)_{[1,i]}^{[1,i]}}{\sigma_i}.
\]
This is obtained from
\eqref{eq:phi det} and
\eqref{eq:PhiL_Z}.

Let us study the matrix
$
\Xi:=
(E-{\rm e}^{a_1}C_A)P^{-1}ZP$.
Because the entries of $C_A$ and $P^{-1}ZP$ are $S_n$-invariant (see Remark~\ref{rem:PinvZP}),
we have
$
 \omega^i \Xi=
(E-{\rm e}^{a_{i+1}}C_A)P^{-1}ZP$.

In order to finish the proof of Proposition~\ref{prop:phi}, it suffices to prove the following.
\begin{Lemma}\label{lem:det Xi}
Let \smash{$M_{(n-i-1)^i}'$} be as defined in Corollary~{\rm\ref{cor:g_rectangle}}.
We have
\begin{equation}
\det \Xi_{[1,i]}^{[1,i]}=\omega^{i+1} \det M_{(n-i-1)^i}'.
\label{eq:det Xi}
\end{equation}
\end{Lemma}
\begin{proof}
Recall $C_A=P^{-1}AP$.
Because $P^{-1}$ is lower unitriangular, we have
\[
\det\Xi^{[1,i]}_{[1,i]}
=
\det((E-{\rm e}^{a_1}A)ZP)_{[1,i]}^{[1,i]}.
\]
As all elements in the $1$st column of $E-{\rm e}^{a_1}A$ are zero and \smash{$(E-{\rm e}^{a_1}A)_{[1,n-1]}^{[2,n]}$} is lower triangular, the determinant \smash{$\det((E-{\rm e}^{a_1}A)ZP)_{[1,i]}^{[1,i]}$} decomposes as
\[
\det(E-{\rm e}^{a_1}A)^{[2,i+1]}_{[1,i]}\det(ZP)_{[2,i+1]}^{[1,i]}
={\rm e}^{ia_1}\det(ZP)_{[2,i+1]}^{[1,i]}.
\]
Therefore, we obtain
\begin{align}
\det\Xi^{[1,i]}_{[1,i]}
&={\rm e}^{ia_1}\det(ZP)_{[2,i+1]}^{[1,i]}
={\rm e}^{ia_1}\det\left(
\begin{matrix}
Z_{[2,i+1]}^{[1,n]} \\\hline
\bigl(P^{-1}\bigr)_{[i+1,n]}^{[1,n]}
\end{matrix}
\right), \label{eq:Ni det}
\end{align}
where for the last equality we used $\det(P)=1$ (cf.\ proof of Corollary~\ref{cor:k rect and tau}).
By row reduction, we can replace \smash{$\bigl(P^{-1}\bigr)_{[i+1,n]}^{[1,n]}$} in~\eqref{eq:Ni det} by
\begin{equation}
\left(
\begin{array}{@{}c|c@{}}
 (-{\rm e}^{-a_1})^{i} &\ast\ \cdots\ \ast \\\hline
 \bm{0}_{n-i-1} & \omega\bigl(P^{-1}\bigr)_{[i+1,n-1]}^{[1,n-1]}
\end{array}
\right).
\label{eq:sweep}
\end{equation}
In fact we sweep out the first column of \smash{$\bigl(P^{-1}\bigr)_{[i+1,n]}^{[1,n]}$} taking $(1,1)$ entry $(-{\rm e}^{-a_1})^i$ as
the pivot. By using the following obvious identity
\[
h_m(x_1,\dots,x_n)-h_{m-1}(x_1,\dots,x_n)x_n=h_{m}(x_1,\dots,x_{n-1})
\] repeatedly, we obtain~\eqref{eq:sweep}.
Thus by expanding along the first column, we have
\[
\det\Xi^{[1,i]}_{[1,i]}=
\det\left(
\begin{matrix}
Z_{[2,i+1]}^{[2,n]} \\\hline
\omega\bigl(P^{-1}\bigr)_{[i+1,n-1]}^{[1,n-1]}
\end{matrix}
\right).
\]
If we compare this with
\[M_{(n-i-1)^i}'=
\omega^{-i} \left(
\begin{matrix}
Z^{[1,n-1]}_{[1,i]}\\
\hline
\bigl(P^{-1}\bigr)^{[1,n-1]}_{[i+1,n-1]}
\end{matrix}
\right),
\]
and use~\eqref{eq:pi_zij},
we obtain~\eqref{eq:det Xi}.
\end{proof}

\begin{proof}
[Proof of Proposition~\ref{prop:phi}]
\begin{align*}
\tilde{\Phi}_n(\fac_i)&=
{\sigma_i}^{-1}{\det}\bigl(\omega^{i} \Xi\bigr)^{[1,i]}_{[1,i]}={\sigma_i}^{-1}{\det}\bigl(\omega^{2i+1} M'_{(n-i-1)^i}\bigr)^{[1,i]}_{[1,i]} \qquad\text{by Lemma~\ref{lem:det Xi}}\\
&={\sigma_i}^{-1}
\omega^{2i+1}\bigl({\co_{(n-i-1)^i}(y)}\cdot{\tilde{g}^{(k)}_{(n-i-1)^i}(y|b)}\bigr)\qquad\text{by Corollary~\ref{cor:g_rectangle}}\\
&={\sigma_i}^{-1}\prod_{l=1}^i\Omega(b_{i+l+1}|y)\cdot\tilde{g}^{(k)}_{(n-i-1)^i}\bigl(y|\omega^{2i+1}b\bigr).\tag*{\qed}
\end{align*}\renewcommand{\qed}{}
\end{proof}

\subsection{Proof of Theorem~\ref{thm:max factor}}
\begin{Lemma}\label{lem:rho check}
Let \smash{$\rho^\lor=\sum_{i=1}^{n-1}\varpi_i^\lor$}.
Suppose $n$ is even and $n=2m$. Then
$\rho^\lor+\varpi_m^\lor\in Q^{\lor}$ and
$
w_\circ t_{-\rho^\lor-\varpi_m^\lor}
=x_{R_m\cup \nu_n}$.
Suppose $n$ is odd. Then $\rho^\lor\in Q^{\lor}$, and
$
w_\circ t_{-\rho^\lor}
=x_{\nu_n}$.
\end{Lemma}
\begin{proof}
It is well known that
\begin{equation}\label{eq:rho check}
\rho^\lor=\frac{1}{2}\sum_{i=1}^{n-1
}i(n-i)\alpha_i^\lor.
\end{equation}
We see that if $n$ is odd then $\rho^\lor\in Q^\lor$
and if $n$ is even then
$\rho^\lor+\varpi_m^\lor\in Q^\lor$.

We work in the extended affine Weyl group $\tilde{W}_G:=\langle \pi\rangle\ltimes \hat{W}_G$
where $\pi^n=\mathrm{id}$ and $\pi s_i=s_{i+1}\pi$ for $i\in \tilde{I}=\Z/n\Z$.
By using \cite[Lemma 4.6]{IIN}, we have
\begin{equation}
\label{eq:w0t}
w_\circ t_{-\rho^\lor}=
w_\circ t_{-\sum_{i=1}^{n-1}\varpi_i^\lor}=
\pi^{-\sum_{i=1}^{n-1}i}x_{\nu_n}.
\end{equation}
Suppose $n=2m$ is even.
Then \smash{$
\pi^{-\sum_{i=1}^{n-1}i}=\pi^{m}$}.
From this, together with
\smash{$t_{-\varpi_m^\lor}=\pi^{-m} x_{R_m}$},
we have
\[
w_\circ t_{-\rho^\lor-\varpi_m^\lor}
=\pi^{m}x_{\nu_n}\pi^{-m} x_{R_m}
=x_{\nu_n\cup R_m},
\]
where, for the last equality, we use \cite[Lem\-ma~2.15]{Tak}.
If $n$ is odd, then
\smash{$
\pi^{-\sum_{i=1}^{n-1}i}=\mathrm{id}$} and the desired equality is nothing but~\eqref{eq:w0t}.
\end{proof}

\begin{proof}[Proof of Theorem~\ref{thm:max factor}]
We first note that $\nu_n$
and $R_m$ are invariant under $\omega_k$.

Suppose $n$ is odd.
From Theorem~\ref{thm:main},
\[
\tilde{\Phi}_n\bigl(Q^{-\rho^\lor}\mathfrak{G}_{w_\circ}^Q\bigr)=\tilde{g}_{\nu_n^{\omega_k}}^{(k)}(y|b)=
\tilde{g}_{\nu_n}^{(k)}(y|b).
\]
From
\eqref{eq:rho check}, it is straightforward to show
\[
\tilde{\Phi}_n\bigl(Q^{\rho^\lor}\bigr)=
\frac{\sigma_{n}^{(n-1)/2}}{\sigma_1\sigma_2\cdots \sigma_{n-1}}.
\]
Comparing these equations with
Proposition~\ref{prop:phi}
and using
$\sigma_n=\Omega(b_1|y)\cdots\Omega(b_n|y)$, we obtain~\eqref{eq:odd}.

Next we consider the case
$n=2m$. From Theorem~\ref{thm:main}, we have
\begin{align*}
\tilde{\Phi}_n
\bigl(Q^{-\rho^\lor-\varpi_m^\lor}
\mathfrak{G}_{w_\circ}^Q\bigr)
&=\tilde{g}^{(k)}_{(\nu_n\cup R_m)^{\omega_k}}(y|b)\qquad \text{by Lemma~\ref{lem:rho check}}\\
&=\tilde{g}_{\nu_n\cup R_m}^{(k)}(y|b)=\tilde{g}^{(k)}_{R_m}(y|b)\tilde{g}^{(k)}_{\nu_n}(y|\omega^m b)\qquad\text{by Theorem~\ref{thm:$k$-rect}.}
\end{align*}
From~\eqref{eq:rho check}, it is straightforward to show
\[
\tilde{\Phi}_n\bigl(Q^{\rho^\lor+\varpi_m^\lor}\bigr)=
\frac{\sigma_{n}^m}{\sigma_m(\sigma_1\sigma_2\cdots \sigma_{n-1})}.
\]
Comparing these equations with
Proposition~\ref{prop:phi}, and
using $\sigma_n=\Omega(b_1|y)\cdots\Omega(b_n|y)$ together with
\smash{$\tilde{g}_{R_m}^{(k)}(y|b)=\sigma_m/\co_{R_m}(y)$}, we obtain~\eqref{eq:even}.
\end{proof}

\appendix

\section{Discrete relativistic Toda lattice}\label{sec:Hirotas_bilinear_form}

Let $L=MN^{-1}$ be the Lax matrix of the relativistic Toda equation~\eqref{eq:rel_Toda}.
Define the new matrix $L^+=N^{-1}M$ by switching the position of $M$ and $N^{-1}$.
By the Gauss decomposition~${L^+=M^+\bigl(N^+\bigr)^{-1}}$, we can define the rational map $(N,M)\mapsto \bigl(N^+,M^+\bigr)$, where $N^+$ and~$M^+$ are matrices of the form
\[
N^+=
\begin{pmatrix}
1& \\
-Q^+_1z^+_1 & 1\\
 & \!\!-Q^+_2z^+_2 & \ddots\\
 & & \ddots & 1\\
 & & & -Q^+_{n-1}z^+_{n-1} &1
\end{pmatrix},\qquad
M^+=
\begin{pmatrix}
z^+_1& -1\\
 & z^+_2 &\ddots\\
 & & \ddots & -1\\
 & & & z^+_n
\end{pmatrix}.
\]

Indeed, this correspondence defines an automorphism over the quotient field of~$R[T][z,Q]$, written as
$
\mathrm{dToda}\colon \CC(z_i,Q_i)\to \CC(z_i,Q_i)$, $ z_i\mapsto z_i^+$, $ Q_i\mapsto Q_i^+$.
This birational map is explicitly written as
\[
Q^+_i=\frac{z_{i}}{z_{i+1}}Q_i, \qquad
z^+_i=\frac{1-Q^+_{i-1}}{1-Q^+_{i}}z_i
, \qquad Q_0:=0,
\]
which is known as \textit{the discrete Relativistic Toda lattice} \cite{Suris}.

As $NL^+=LN$, the two matrices $L$ and $L^+$ are similar.
Following the same argument as~\eqref{eq:Gauss_decomp} and~\eqref{eq:AZ_to_L}, we have the pair of matrices $U^+,R^+$ satisfying
\[
P^{-1}Z^+AP=\bigl(U^+\bigr)^{-1}R^+\qquad \mbox{and} \qquad L^+=U^+C_A\bigl(U^+\bigr)^{-1}=R^+C_A(R^+)^{-1}.
\]
Comparing them with $L^+=N^{-1}LN=M^{-1}LM$, we obtain the equations
$
U^{+}=N^{-1}U$, $ R^+=c \cdot M^{-1}R
$
with some nonzero constant $c$.
From these expressions, we can show $P^{-1}Z^+AP=c\cdot P^{-1}ZP$, and then,
$
Z^+=cZA^{-1}$.
The correspondence $Z\mapsto Z^+$ defines the same map as $\sigma$ (see~\eqref{eq:def_of_sigma}) on the centralizer space $\tilde{\mathscr{Z}}$.
Symbolically, we have $\sigma=\mathrm{dToda}$.

Let
$\tau^+_i:=\sigma_i=\sigma(\tau_i)$
and $\tau^-_i:=\sigma^{-1}(\tau_i)$ be the tau-functions at the next and previous time step, respectively.
Substituting~\eqref{eq:Q_to_tau} and~\eqref{eq:z_to_tau} to
\[
z_i^+=\frac{1-Q^+_{i-1}}{1-Q^+_{i}}z_i,
\] we obtain
\[
\frac{\tau_i^2-\tau_{i+1}\tau_{i-1}}{\tau_i^+\tau_i^-}=
\frac{\tau_{i-1}^2-\tau_{i}\tau_{i-2}}{\tau_{i-1}^+\tau_{i-1}^-}=\cdots
=
\frac{\tau_{1}^2-\tau_{2}\tau_{0}}{\tau_{1}^+\tau_{1}^-}=1.
\]
Finally, we obtain
$
\tau_i^2-\tau_{i+1}\tau_{i-1}=\tau_i^+\tau_i^-$,
which is known as \textit{Hirota's bilinear form of the discrete Toda equation}.

\section{Proof of (\ref{eq:polynomialF})}
\label{sec:F}
Note that for any square matrix $A$ of size $n$, we have
\[
\det(\zeta E-A)=
\sum_{i=0}^n(-1)^i\cdot \zeta^{n-i}
\sum_{
\substack{
J\subset [1,n]\\
|J|=i
}
}
\det A_J^J.
\]
Therefore, it suffices to show
\[
\det L_J^J=
\prod_{j\in J,j+1\notin J}(1-Q_j)
\prod_{j\in J}z_j.
\]
Because the matrix $L^{[a,b]}_{[a,b]}$ decomposes as
\[
\begin{pmatrix}
z_a & -1 && & \\
 & z_{a+1}&\ddots& &\\
 & &\ddots& -1 &\\
 & &&\!z_{b-1}&-1\\
 &&&&\!\!\!\!\!(1-Q_b)z_b
\end{pmatrix}
\begin{pmatrix}
1 & & && & \\
-Q_az_a & 1 && & \\
 &\!\!\!\!\!\!\!-Q_{a+1}z_{a+1} &\ddots& &\\
 & & \ddots &1&\\
 &&&\!\!\!\!\!\!\!\!\!-Q_{b-1}z_{b-1}&1
\end{pmatrix}^{-1},
\]
we have
\[
\det L^{[a,b]}_{[a,b]}=
(1-Q_b)z_a\cdots z_b
\qquad
\text{for}\quad 1\le a\le b\le n.
\]
If we decompose
$J$ as
$J=J_1\sqcup \cdots\sqcup J_r$ in such a~way that
$J_i=[a_i,b_i]$
with $a_{i+1}-b_i\ge 2 $, then
$L_J^J$ is blockwise lower triangular with respect to
$J=J_1\sqcup \cdots\sqcup J_r$. So we have
$\det L_J^J
=\prod_{i=1}^r \det L_{J_i}^{J_i}$
which concludes~\eqref{eq:polynomialF}.

\section[Affine K-nil-Hecke action on quantum K-ring]{Affine $\boldsymbol{K}$-nil-Hecke action on quantum $\boldsymbol{ K}$-ring}\label{sec:nil-Hecke_on_QK}
The aim of this appendix is to give a~proof of Proposition~\ref{prop:Quantum_D(1)}. We use the theory of
semi-infinite flag manifolds.

Let $G$ be a~simple, simply-connected complex algebraic group. Let $B$ be a~Borel subgroup of $G$ and $T$ a~maximal torus contained in $B$.
Let $Q$, $P$ be the root lattice and weight lattice, respectively.
Let $W_G=\langle s_i \mid i\in I\rangle$ be the Weyl group of $G$, and $\hat{W}_G=W_G\ltimes Q^\lor$ the affine Weyl group with $\hat{W}_G=\big\langle s_i \mid i\in \Iaf\big\rangle$.
The \emph{semi-infinite flag manifold} associated with $G$
is the (reduced) ind-scheme of ind-infinite type \smash{$\mathbf{Q}_G^{\mathrm{rat}}= G(\mathbb{F})/(T(\C)N(\mathbb{F}))$} (see~\cite{K:Fro,KNS}), where $\mathbb{F}=\C(\!(t)\!)$ and~$N$ is the unipotent radical of $B$. For each $x\in \hat{W}_G$, there is the corresponding Schubert variety~\smash{$\mathbf{Q}_G(x)\subset \mathbf{Q}_G^{\mathrm{rat}}$}. Let $\mathbf{Q}_G=\mathbf{Q}_G(e)$, where $e$ is the identity.
Let $Q^{\lor,+}$ denote the positive part of the coroot lattice.
The $T$-equivariant $K$-group $K_T(\mathbf{Q}_G)$ is, as an $R(T)$-module, the direct product \smash{$\prod_{x\in \hat{W}_{G}^{+}}R(T)\O_{\mathbf{Q}_G(x)}$}, where~${\hat{W}_{G}^{+}=
\big\{wt_\xi\mid w\in W_G,\, \xi\in Q^{\lor,+}\big\}}$.
Kato~\cite{Kato} established an $R(T)$-module isomorphism from
$QK_T(G/B)$ to $K_T(\mathbf{Q}_G)$.
Explicitly,
${\rm e}^\mu \O^wQ^\xi$ in $QK_T(G/B)$, with $\mu\in P$, $w\in W_G$, $\xi \in Q^\lor$,
corresponds to
${\rm e}^{-\mu}\O_{\mathbf{Q}_G(wt_\xi)}$
in $K_T(\mathbf{Q}_G)$; here we follow the convention of \cite[Section~6.1]{LNS:Se} and \cite[Section~5]{MNS1}.

The \emph{semi-infinite length} of $x=wt_\xi\in \hat{W}_G$
is defined by
\smash{$\ell^\sinf(w)
=\ell(w)+2\langle \rho,\xi\rangle$}, where $\rho$ is~the half sum of the positive roots.
Let \smash{$\le_\sinf$} denote the \emph{semi-infinite Bruhat order} (see \cite[Section~2.4]{INS}).
We note that
\[
s_i x>_\sinf x\Longleftrightarrow \ell^\sinf(s_ix)=\ell^\sinf(x)+1, \qquad i\in \tilde{I}
\] (see \cite[Lemma~4.1.2]{INS}).
Let $\hat{W}_G^0$ denote the minimal-length coset representatives for the coset space $\hat{W}_G/W_G$
with respect to the ordinary length function.

The following fact is due to Peterson~\cite{Pet}.
\begin{Lemma}\label{lem:sinf}
For \smash{$x\in \hat{W}_G^0$}, we have
\smash{$\ell(x)=-\ell^{\sinf}(x)$}.
\end{Lemma}
\begin{proof}
Let $x=wt_\lambda\in \hat{W}_G^0$ $\bigl(w\in W_G, \, \lambda\in Q^\vee\bigr)$. Note that $\lambda\in Q^\vee$ is anti-dominant according to \cite[Lemma 3.3]{LS:Acta}.
We have
\begin{align*}
 \ell(x)&=\ell(t_\lambda)-\ell(w)=\langle e\cdot \lambda,-2\rho\rangle-\ell(w)=-\ell^\sinf(x),
\end{align*}
where the first equality uses \cite[Lemma 3.3]{LS:Acta}, and the second equality applies \cite[Lemma 3.2]{LS:Acta} with $w=e$.
\end{proof}

Motivated by the nil-DAHA action on $K_T\bigl(\mathbf{Q}_G^{\mathrm{rat}}\bigr)$
\cite[Theorem 6.5]{KNS} (see also \cite[Section~3.1]{KoNOS}, \cite[Section~2.6]{Orr}),
we define
an endomorphism of
$QK_T(G/B)_Q$,
\[
D_0^Q=T_\theta+
Q^{-\theta^\vee}
\mathscr{O}^{s_\theta}\cdot
s_\theta,
\]
where~$\theta$ is the highest root and $T_\theta$ is defined in~\eqref{eq:Ttheta}.
Then $D_i$ $ \big(i\in \tilde{I}\big)$ satisfy the braid relations, and hence $D_x^Q$ is defined for $x\in \hat{W}_G$.

The next result is transported from
the corresponding one for semi-infinite flag manifolds
(see \cite[Section~2.6]{Orr} and \cite[Section~3.1]{KoNOS}).
\begin{Proposition}\label{prop:left_Dem}
Let $w\in W$. For $i\in I$, we have
\begin{align}
\label{eq:DiQGro}
D_i^Q(\mathscr{O}^w)
&=\begin{cases}
\mathscr{O}^{s_iw} & \text{if}\ s_iw<w,\\
 \mathscr{O}^{w} & \text{if}\ s_iw>w.
\end{cases} \end{align}
Moreover, we have
\begin{align}
\label{eq:D0QGro}
D_0^Q(\mathscr{O}^w)
&=\begin{cases}
 Q^{-w^{-1}(\theta^\vee)}\mathscr{O}^{s_\theta w} & \text{if}\ s_\theta w>w,\\
 \mathscr{O}^{w} & \text{if}\ s_\theta w<w.
\end{cases}
\end{align}
\end{Proposition}
\begin{proof}
Let $x=wt_\xi$
with $w\in W_G$
and $\xi\in Q^{\lor}$.
We use the following fact (see \cite[Appendix~A]{INS})
\begin{align*}
\text{for}\quad i\in I,\qquad s_i x<_\sinf x   \ \Longleftrightarrow \ s_i w< w,\qquad
 s_0 x<_\sinf x \ \Longleftrightarrow \ s_\theta w>w.
\end{align*}
By the isomorphism
$QK_T(G/B)\cong K_T(\mathbf{Q}_G)$ described above,
\eqref{eq:DiQGro} and~\eqref{eq:D0QGro} are obtained from
\cite[equation~(2.28)]{Orr} or \cite[equation~(3.19)]{KoNOS}.
Note that $s_0w=s_\theta w t_{-w^{-1}(\theta^\lor)}$
for $w\in W_G$.
\end{proof}

 \begin{Corollary} \label{cor:BF}
Let $x=wt_\xi\in \hat{W}_G^0$. Then
$
D_x^Q(1)
=Q^\xi \mathscr{O}^{w}$.
\end{Corollary}
\begin{proof}
Let $x=s_{i_1}\cdots s_{i_l}$ be a~reduced expression.
Then we have a~saturated decreasing chain of elements in
$\hat{W}_G^0$
\[
x=s_{i_1}\cdots s_{i_l}
>s_{i_2}\cdots s_{i_l}
>
\cdots
>s_{i_l}>e.
\]
By Lemma~\ref{lem:sinf}, we obtain a~saturated increasing chain of elements
with respect to $<_\sinf$
\[
x=s_{i_1}\cdots s_{i_l}
<_\sinf s_{i_2}\cdots s_{i_l}
<_\sinf
\cdots
<_\sinf s_{i_l}<_\sinf e.
\]
Thus from Proposition~\ref{prop:left_Dem}, we deduce the corollary.
\end{proof}
\begin{proof}[Proof of Proposition~\ref{prop:Quantum_D(1)}]
Apply Corollary~\ref{cor:BF} to the $T$-equivariant
quantum $K$-ring of \linebreak
$\SL_n(\C)$ presented in~\cite{MNS1} and use the fact that
$\Gro{w}(z|\eta)$ represents $\O^w$.
\end{proof}

\subsection*{Acknowledgements}The authors are grateful to the anonymous referees for valuable comments and suggestions, which helped improve the manuscript.
We are grateful for the fruitful discussions with Mark Shimozono, Toshiaki Maeno, Takafumi Kouno, and Daisuke Sagaki.
T.I.\ was partly supported by JSPS Grants-in-Aid for Scientific Research 23H01075, 22K03239, 20K03571 and 20H00119.
S.I.\ was partly supported by JSPS Grants-in-Aid for Scientific Research 22K03239 and 23K03056.
S.N.\ was partly supported by JSPS Grant-in-Aid for Scientific Research 21K03198.

\pdfbookmark[1]{References}{ref}
\LastPageEnding

\end{document}